\documentclass[12pt,reqno]{amsart}
\theoremstyle{plain}
\newtheorem{theorem}{Theorem}[section]
\newtheorem{corollary}[theorem]{Corollary}
\newtheorem{example}[theorem]{Example}

\newtheorem{proposition}[theorem]{Proposition}
\newtheorem{lemma}[theorem]{Lemma}

\theoremstyle{definition}
\newtheorem{definition}[theorem]{Definition}

\theoremstyle{remark}
\newtheorem{remark}[theorem]{Remark}
\newtheorem{problem}[theorem]{Problem}

\begin{document}

\title{Frames in Hilbert C*-modules and C*-algebras}
\author[M.~Frank]{Michael Frank}
\address{Universit\"at Leipzig, Math.~Institut, D-04109 Leipzig, F.R.~Germany}
\email{frank@mathematik.uni-leipzig.de}
\thanks{Both authors were supported by grants from the NSF}
\author[D.~R.~Larson]{David R.~Larson}
\address{Dept.~Mathematics, Texas A{\&}M Univ., College Station, TX 77843, U.S.A.}
\email{larson@math.tamu.edu}
%%%%%%%%%%%%%%%%%%%%%%%%%%%%%%%%%%%%%%%%%%%%%%%%%%%%%%%%%%%%%%%%%%%%%%%%%%%%%
\begin{abstract}
We present a general approach to a module frame theory in C*-algebras and
Hilbert C*-modules. The investigations rely on the ideas of geometric dilation
to standard Hilbert C*-modules over unital C*-algebras that possess orthonormal
Hilbert bases, of reconstruction of the frames by projections and by other
bounded module operators with suitable ranges. We obtain frame representation
and decomposition theorems, as well as similarity and equivalence results.
Hilbert space frames and quasi-bases for conditional expectations of finite
index on C*-algebras appear as special cases. Using a canonical
categorical equivalence of Hilbert C*-modules over commutative C*-algebras
and (F)Hilbert bundles the results find a reinterpretation for frames in
vector and (F)Hilbert bundles.
\end{abstract}
%%%%%%%%%%%%%%%%%%%%%%%%%%%%%%%%%%%%%%%%%%%%%%%%%%%%%%%%%%%%%%%%%%%%%%%%%%%%%
\maketitle

%%%%%%%%%%%%%%%%%%%%%%%%%%%%%%%%%%%%%%%%%%%%%%%%%%%%%%%%%%%%%%%%%%%%%%%%%%%%%

The purpose of this paper is to extend the theory of frames known for
(separable) Hilbert spaces to similar sets in C*-algebras and (finitely and
countably generated) Hilbert C*-modules. The concept 'frame' may generalize
the concept 'Hilbert basis' for Hilbert C*-modules in a very efficient way
circumventing the ambiguous condition of 'C*-linear independence' and emphasizing
geometrical dilation results and operator properties. This idea is natural in
this context because, while such a module may fail to have any reasonable type
of basis, it turns out that countably generated Hilbert C*-modules over unital
C*-algebras always have an abundance of frames of the strongest (and simplest)
type. The considerations follow the line of the geometrical and
operator-theoretical approach worked out by Deguang Han and David R.~Larson
\cite{HanLarson} in the main. They include the standard Hilbert space case in
full as a special case, see also
\cite{DLS1:97,DLS2:98,GH:98,HL:98,ILP:98,WUTAM}. However, proofs that
generalize from the Hilbert space case, when attainable, are usually
considerably more difficult for the module case for reasons that do not occur
in the simpler Hilbert space case. For example, Riesz bases of Hilbert spaces
with frame bounds equal to one are automatically orthonormal bases, a straight
consequence of the frame definition. A similar statement for standard Riesz
bases of certain Hilbert C*-modules still holds, but the proof of the statement
requires incomparably more efforts to be established, see Corollary
\ref{on-basis}.
Generally speaking, the known results and obstacles of Hilbert C*-module
theory in comparison to Hilbert space and ideal theory would rather suggest
to expect a number of counterexamples and diversifications of situations
that could appear investigating classes of Hilbert C*-modules and of
C*-algebras of coefficients beyond the Hilbert space situation. Surprisingly,
almost the entire theory can be shown to survive these significant changes.
For complementary results to those explained in the present paper we refer to
\cite{FL99,FL00}.

This paper began with a talk the second author gave on the content of
\cite{HanLarson} at the Joint COAS-GPOTS symposium at Kingston, Ontario,
in May 1997. After the talk the first author suggested that many of the ideas
concerning frames in the Hilbert space situation may have natural counterparts
in Hilbert C*-module theory. He proposed that we consider joint work attempting
to use \cite{HanLarson} as a 'blueprint' for ideas. After we got deeply into
the project we discovered that frames and related ideas had in fact been used
by others implicitly and explicitly in the C*-literature (although the term
'frame' had not been applied, the connection with engineering literature
had not been realized, and the constructions and ideas had not been
systematically explored by the authors). At the other side operator-valued
inner products appeared as arguments in proofs of wavelet theory publications
without any reference to Hilbert C*-module theory. (Detailed references will be
given below at the end of the introduction.)

\smallskip
The areas of applications indicate a large potential of problems for the
investigation of which our results could be applied. For example, an
interpretation of our results in terms of noncommutative geometry leads to
frames in vector bundles and (F)Hilbert bundles, \cite{Serre,Swan,DuGi}.
The decision was made
to publish the core results of our work in a way that should bring them to
the attention of an audience beyond researchers working in the field of
operator theory and operator algebras. So some of the explanations in the
following sections may contain some more details than specialists may need 
to understand the presented theory.

\smallskip
By the commonly used definition of a (countable) frame in a (separable)
Hilbert space a set $\{ x_i : i \in \mathbb J \} \subset H$ is said to
be a frame of the Hilbert space $H$ if there exist two constants $C$, $D > 0$
such that the inequality
\[
  C \cdot \| x \|^2 \leq {\sum}_i | \langle x,x_i \rangle |^2
  \leq D \cdot \|x\|^2
\]
holds for every $x \in H$. To generalize this definition to the situation of
Hilbert C*-modules we have to rephrase the inequality in a suitable way.
Therefore, frames of Hilbert $A$-modules $\{ {\mathcal H}, \langle .,. \rangle
\}$ over unital C*-algebras $A$ are sets of elements $\{ x_i : i \in {\mathbb
J} \} \subset \mathcal H$ for which there exist constants $C$, $D > 0$ such
that the inequality
  \begin{equation}  \label{eq00}
   C \cdot \langle x,x \rangle \leq {\sum}_i \langle x,x_i
   \rangle   \langle x_i,x \rangle \leq D \cdot \langle x,x \rangle
  \end{equation}
is satisfied for every $x \in \mathcal H$. An additional restriction to
the sum in the middle of the inequality (\ref{eq00}) to converge in norm for
every $x \in \mathcal H$ guarantees the existence and the adjointability of
the frame transform $\theta: {\mathcal H} \to l_2(A)$ and the orthogonal
comparability of its image inside $l_2(A)$, facts that are crucial and
unexpected in the generality they hold. The restriction to countable frames
is of minor technical importance, whereas the restriction to unital C*-algebras
of coefficients refers to the fact that approximative identities of
non-unital C*-algebras do not serve as approximative identities of their
unitizations. The investigation of arbitrary frames with weakly converging sums
in the middle of (\ref{eq00}) requires Banach C*-module and
operator module techniques and has to be postponed. Some remarks on this
problem are added in section eight of the present paper. We point out that frames
exist in abundance in finitely or countably generated Hilbert C*-modules over
unital C*-algebras $A$ as well as in the C*-algebras itself, see Example
\ref{ex4}. This fact allows to rely on standard decompositions for elements
of Hilbert C*-modules despite of the general absence of orthogonal and
orthonormal Riesz bases in them, cf.~Example \ref{ex214}.

\smallskip
From the point of view of applied frame theory the advantage of the generalized
setting of Hilbert C*-modules may consist in the additional degree of
freedom coming from the C*-algebra $A$ of coefficients and its special inner
structure, together with the handling of the basic features of the generalized
theory in almost the same manner as for Hilbert spaces. For example, for
commutative C*-algebras $A = {\rm C}(X)$ over compact Hausdorff spaces $X$,
continuous (in some sense) fields of frames over $X$ in the Hilbert space
$H$ could be considered using the geometric analogues of Hilbert
${\rm C}(X)$-modules - the vector bundles or (F)-Hilbert bundles with base
space $X$. An appropriate choice of the compact base space of the bundles
allows the description of parameterized and continuously varying families of
classical frames in a given Hilbert space.

\smallskip
The content of the present paper is structured as follows: Section 1 contains
the preliminary facts about Hilbert C*-module theory needed to explain our
concept. Section 2 covers the definition of the different types of frames
in C*-algebras and Hilbert C*-modules and explains some of their basic
properties. Section 3 is devoted to a collection of re\-presentative examples
showing the phenomena that have to be taken into account for a generalization
of the theory away from Hilbert spaces to Hilbert C*-modules. The existence of
the frame transform $\theta$, its properties and the reconstruction formula
for standard normalized tight frames are proved in section 4 giving the key
to a successful generalization process. In particular, standard normalized tight
frames are shown to be sets of generators for the corresponding Hilbert
C*-modules. In section 5 geometrical dilation results and similarity problems
of frames are investigated and results are obtained covering the general
situation. The existence and the properties of canonical and alternate dual
frames is the goal of section 6. As a consequence a reconstruction formula
for standard frames is established. The last section contains a classification
result showing the strength of the similarity concept of frames.
Some final remarks complete our investigations.

In the present paper some results have been obtained for the theory of Hilbert
C*-modules which are partially new to the literature and which use our frame
technique in their proofs, see the Propositions \ref{prop-magic}, \ref{prop-HS} and Theorem
\ref{prop-afgHm}. In particular, we prove that every set of algebraic generators
of an algebraically finitely generated Hilbert C*-module is automatically
a module frame. We give a new short proof that any finitely generated
Hilbert C*-module is projective. Beside this, a new characterization of
Hilbert-Schmidt operators on Hilbert spaces allows to extend this concept to
certain classes of Hilbert C*-modules over commutative C*-algebras.

\smallskip
At this place we want to give more detailed references to the literature to
appreciate ideas and work related to our results that have been published by
other researchers. Most of the publications listed below were not known to us
at the time we worked out modular frame theory in 1997-1998. Some of the
mentioned articles have been written very recently.

We make use of G.~G.~Kasparov's Stabilization Theorem (\cite[Th.~1]{Kas}) in an
essential way. However, far not every set of generators of countably generated
Hilbert C*-modules admits the frame property, even in the particular situation
of separable Hilbert spaces. Our aim is to divide out this special class of
generating sets and to characterize them as powerful structures in
countably generated Hilbert C*-modules that are capable to play the role
bases play for Hilbert spaces.

Another source of inspiration has been the inner structure of self-dual
Hilbert W*-modules described by W.~L.~Paschke in \cite{Pa1} in 1973.
Rephrasing his description in the context of frames it reads as the proof of
general existence of orthogonal normalized tight frames $\{ x_j : j \in
\mathbb J \}$ for self-dual Hilbert W*-modules, where additionally the values
$\{ \langle x_j,x_j \rangle : j \in \mathbb J \}$ are projections. This point
of view has been already realized by Y.~Denizeau and J.~F.~Havet in \cite{DH}
in 1994 as pointed out to us by the referee. They went one step further taking
a topologically weak reconstruction formula for normalized tight frames as
a corner stone to characterize the concept of `quasi-bases' for Hilbert
W*-modules. The special frames appearing from W.~L.~Paschke's result are
called `orthogonal bases' by these authors. The two concepts have been
investigated by them to the extend of tensor product properties of quasi-bases
for C*-correspondences of W*-algebras,
cf.~\cite[Thm. 1.2.5, Cor.~1.2.6, Lemma 2.1.5]{DH}.
A systematic investigation of the concept of quasi-bases has not been provided
at that place. While these results are surely interesting from the point of
view of operator theory they are only of limited use for wavelet theory. For
our opinion the main reason is the necessity of a number of weak completion
processes to switch from basic Hilbert space contexts to suitable self-dual
Hilbert W*-module contexts. On this way too much structural information
gets lost or hidden, in general.

Looking back into the literature for Y.~Denizeau's and J.-F.~Havet's
motivation to introduce quasi-bases at a rather general level, the concept of
`quasi-bases' can be found to be worked out for the description of
algebraically characterizable conditional expectations of finite index on
C*-algebras by Y.~Watatani in 1990, \cite{Wata}. At that place quasi-bases are
a special example of module frames in Hilbert C*-modules (more precisely, a
pair consisting of a frame and a dual frame). For normal conditional expectations
of finite index on W*-algebras generalized module frames like Pimsner-Popa
bases have been considered earlier by M.~Pimsner and S.~Popa \cite{PP},
by M.~Baillet, Y.~Denizeau and J.-F.~Havet \cite{BDH,DH}, and by E.~Kirchberg
and the author \cite{KiFr98}, among others (cf.~\cite{Pa1,Fr1,Bl2} for
technical background information). Recently, M.~Izumi proved the general existence
of module frames for Hilbert C*-modules that arise from simple C*-algebras
by a conditional expectation of finite index onto one of their C*-subalgebras,
cf.~\cite{Izumi}.
We discovered the use of standard frames in one place of E.~C.~Lance's lecture
notes \cite{Lance2} where he used this kind of sequences in one reasoning
on page 66, without any investigation of the concept itself.
In Hilbert C*-module theory and its applications special generating sequences
have been used to investigate a large class of generalized Cuntz-Krieger-Pimsner
C*-algebras. These C*-algebras arise from Hilbert C*-bimodules in categorical
contexts in the way of making use of existing canonical representations of
elements, \cite[p.~266]{DPZ}, \cite[\S 2]{KPW}. The exploited sequences of
elements of the Hilbert C*-modules under consideration have been called
'bases'. They admit the key frame properties. The authors make use of a
reconstruction formula for bases of that kind, but without any explicit
statement.

We have learned by a communication of M.~A.~Rieffel that the
idea to use finitely generated projective C*-modules over commutative
C*-algebras for the investigation of multiresolution analysis wavelets was
introduced by him in a talk given at the Joint Mathematics Meeting at San
Diego in January 1997, \cite{Rie2}. He has considered module frames generated
by images of a frame in a certain projective C*-submodule and canoni\-cal
representations of elements related to them.
P.~J.~Wood pointed out in \cite[p.~10]{Wood1} that algebra-valued inner
products have been used before by C.~de Boor, R.~DeVore and A.~Ron in 1992,
\cite{BDR}, and by A.~Fischer in 1997, \cite{Fi}. In fact, $L^1$-spaces serve
as target spaces. They used these structures in proofs treating vanishing
moments and approximation properties of wavelets. However, the concept of a
$*$-algebra-valued inner product has not been introduced by these authors.
Similar constructions have been exploited to examine Sobolev smoothness
properties of wavelets, see L.~M.~Villemoes in \cite{Vi} (1992). 

\smallskip
While the present paper has been circulating as a preprint the
ideas and results contained therein have been successfully applied to solve
problems in both operator and wavelet theory.
We know about forthcoming publications by I.~Raeburn and S.~Thompson
\cite{RaeTho} who proved a generalized version of Kasparov's Stabilization
Theorem for a kind of countably generated Hilbert C*-modules over
non-$\sigma$-unital C*-algebras, where the countable sets of generators
consists of multipliers of the module. They generalize our concept of
frames to the situation of certain generating sets consisting of multipliers
of Hilbert C*-modules.
Following the ideas by M.~A.~Rieffel explained in \cite{Rie2} M.~Coco and
M.~C.~Lammers \cite{CoLa} described a W*-algebra and a related self-dual
Hilbert W*-module derived from the analysis of Gabor frames. They showed how
to apply these structures to solve some problems of Gabor analysis. At the
same time P.~J.~Wood analyzed the mentioned ideas by M.~A.~Rieffel in a general
framework of group C*-algebras. Using module frame techniques of Hilbert
C*-module theory he studied the dimension function of wavelets and classified
wavelets by methods derived from C*-algebraic K-theory, see \cite{Wood1,Wood2}.
Motivated by investigations on Hilbert H*-modules D.~Baki\'c and B.~Gulja{\v{s}}
introduced the concept of a `basis' of Hilbert C*-modules over C*-algebras of
compact operators explicitly (i.e.~the concept of normalized tight frames
which are Riesz bases) in 2001, cf.~\cite[Th.~2]{BG}.

\medskip
The authors would like to thank D.~P.~Blecher, P.~G.~Casazza, D.~Kucerovsky,
M.~C.~Lam\-mers, V.~M.~Manuilov, A.~S.~Mishchenko, V.~I.~Paulsen, M.~A.~Rieffel
and E.~V.~Troitsky for helpful conversations and remarks on the subject, as
well as I.~Raeburn, S.~Thompson and P.~J.~Wood for sending copies of their
preprints to us.
They are very grateful to the referee for his conscientious reading of
the presented manuscript and for the valuable remarks and hints in his
review that lead to a substantial improvement of the explanations given below.
The first author is indebted to V.~I.~Paulsen and D.~P.~Blecher for their
invitation to work at the University of Houston in 1998 and for financial
support.

%%%%%%%%%%%%%%%%%%%%%%%%%%%%%%%%%%%%%%%%%%%%%%%%%%%%%%%%%%%%%%%%%%%%%%%%%%%%%%
\section{Preliminaries}

\bigskip
The theory of Hilbert C*-modules generalizes the theory of Hilbert spaces,
of one-sided norm-closed ideals of C*-algebras, of (locally trivial) vector
bundles over compact base spaces and of their noncommutative counterparts -
the projective C*-modules over unital C*-algebras, among others (see
\cite{Lance2,Wegge-Olsen}). Because of the complexity of the theory and
because of the different research fields interested readers of our
considerations may come from we have felt the necessity to give detailed
explanations in places. We apologize to researchers familiar with the basics
of Hilbert C*-module theory for details which may be skipped by more
experienced readers.

\smallskip
Let $A$ be a C*-algebra. A {\it pre-Hilbert $A$-module} is a linear space and
algebraic (left) $A$-module $\mathcal H$ together with an {\it $A$-valued inner
product} $\langle .,. \rangle : {\mathcal H} \times {\mathcal H} \to A$ that
possesses the following properties:

(i) $\,\,\, \langle x,x \rangle \geq 0$ for any $x \in \mathcal H$.

(ii) $\,\, \langle x,x \rangle = 0$ if and only if $x = 0$.

(iii) $\, \langle x,y \rangle = \langle y,x \rangle^*$ for any $x,y
     \in \mathcal H$.

(iv) $\,\, \langle ax+by,z \rangle = a \langle x,z \rangle + b \langle y,z
      \rangle$ for any $a,b \in A$, $x,y,z \in \mathcal H$.

To circumvent complications with linearity of the $A$-valued inner product with
respect to imaginary complex numbers we assume that the linear operations
of $A$ and $\mathcal H$ are comparable, i.e.~$\lambda (ax)=(\lambda a)x=
a(\lambda x)$ for every $\lambda \in {\mathbb C}$, $a \in A$ and $x \in
\mathcal H$. The map $x \in {\mathcal H} \to \| x \| = \| \langle x,x \rangle
\|_A^{1/2} \in {\mathbb R}^+$ defines a norm on $\mathcal H$. Throughout the
present paper we suppose that $\mathcal H$ is complete with respect to that
norm. So $\mathcal H$ becomes the structure of a Banach $A$-module. We refer
to the pairing $\{ {\mathcal H}, \langle .,. \rangle \}$ as to a {\it Hilbert
$A$-module}. Two Hilbert $A$-modules $\{ {\mathcal H}, \langle .,.
\rangle_{{\mathcal H}} \}$ and $\{ {\mathcal K}, \langle .,. \rangle_{{\mathcal
K}} \}$ are {\it unitarily isomorphic} if there exists a bijective bounded
$A$-linear mapping $T: \mathcal H \to \mathcal K$ such that $\langle x,y
\rangle_{\mathcal H} = \langle T(x),T(y) \rangle_{\mathcal K}$ for $x,y \in
\mathcal H$.

\smallskip
If two Hilbert $A$-modules $\{ {\mathcal H}, \langle .,. \rangle_{{\mathcal H}}
\}$ and $\{ {\mathcal K}, \langle .,. \rangle_{{\mathcal K}} \}$ over a
C*-algebra $A$ are given we define their {\it direct sum} ${\mathcal H} \oplus
{\mathcal K}$ as the set of all ordered pairs $\{ (h,k) \, : \, h \in {\mathcal
H}, \, k \in {\mathcal K} \}$ equipped with coordinate-wise operations and with
the $A$-valued inner product $\langle .,. \rangle_{{\mathcal H}} + \langle .,.
\rangle_{{\mathcal K}}$.

\medskip
In the special case of $A$ being the field of complex numbers $\mathbb C$ the
definition above reproduces the definition of Hilbert spaces. However, by far
not all theorems of Hilbert space theory can be simply generalized to the
situation of Hilbert C*-modules. To give an instructive example consider the
C*-algebra $A$ of all bounded linear operators $B(H)$ on a separable Hilbert
space $H=l_2$ together with its two-sided norm-closed ideal $I=K(H)$ of all compact
operators on $H$. The C*-algebra $A$ equipped with the $A$-valued inner
product $\langle .,. \rangle$ defined by the formula $\langle a,b \rangle_A =
ab^*$ becomes a Hilbert $A$-module over itself. The restriction of this
$A$-valued inner product to the ideal $I$ turns $I$ into a Hilbert $A$-module,
too. So we can form the new Hilbert $A$-module ${\mathcal H} = A \oplus I$ as
defined in the previous paragraph. Let us consider some properties of
$\mathcal H$.

First of all, the analogue of the Riesz representation theorem for bounded
($A$-)linear mappings $r: {\mathcal H} \to A$ is not valid for $\mathcal H$.
For example, the mapping $r((a,i)) = a + i$ ($a \in A$, $i \in I$) cannot
be realized by applying the $A$-valued inner product to $\mathcal H$ with one
fixed entry of $\mathcal H$ in its second place since the necessary entry
$(1_A,1_A)$ does not belong to $\mathcal H$.
Secondly, the bounded $A$-linear operator $T$ on $\mathcal H$ defined by the
rule $T: (a,i) \to (i,0_A)$ ($a \in A$, $i \in I$) does not have an adjoint
operator $T^*$ in the usual sense since the image of the formally defined
adjoint operator $T^*$ is not completely contained in $\mathcal H$.
Furthermore, the Hilbert $A$-submodule $I$ of the Hilbert $A$-module $A$ is
not a direct summand, neither an orthogonal nor a topological one.
Considering the Hilbert $A$-submodule ${\mathcal K} \subseteq \mathcal H$
defined as the set ${\mathcal K} = \{ (i,i) : i \in I \}$ with induced from
$\mathcal H$ operations and $A$-valued inner product we obtain the coincidence
of $\mathcal K$ with its biorthogonal complement inside $\mathcal H$. However,
even in this situation $\mathcal K$ is not an orthogonal summand of $\mathcal H$,
but only a topological summand with complement $\{ (a,0_A) : a \in A \}$.

So the reader should be aware that every formally generalized formulation of
Hilbert space theorems has to be checked for any larger class of Hilbert
C*-modules carefully and in each case separately. To provide a collection of
facts from Hilbert C*-module theory used in forthcoming sections the remaining
part of the present section is devoted to a short guideline into parts of the
theory.

\medskip
Let $\mathbb  J$ be a countable set of indices. In case we need a (partial)
ordering on $\mathbb J$ we may choose to identify $\mathbb J$ with the set of
integers $\mathbb N$ or with other countable, partially ordered sets.
A subset $\{ x_j : j \in {\mathbb J} \}$ of a Hilbert $A$-module $\{
{\mathcal H}, \langle .,. \rangle \}$ is a {\it set of generators of
$\mathcal H$} (as a Banach $A$-module) if the $A$-linear hull of $\{ x_j :
j \in {\mathbb J} \}$ is norm-dense in $\mathcal H$. The subset $\{ x_j :
j \in {\mathbb J} \}$ is {\it orthogonal} if $\langle x_i,x_j \rangle = 0$
for all $i,j \in {\mathbb J}$ with $i \not= j$. A set of generators $\{ x_j :
j \in {\mathbb J} \}$ of $\mathcal H$ is a {\it Hilbert basis of $\mathcal H$}
if $\,$ (i) $\, A$-linear combinations $\sum_{j \in S} a_jx_j$ with
coefficients $\{ a_j \}$ in $A$ and $S \subseteq {\mathbb J}$ are equal to
zero if and only if in particular every summand $a_jx_j$ equals to zero for
$j \in S$, and $\,$ (ii) $\, \| x_j \| =1$ for every $j \in {\mathbb J}$.
This definition is consistent since every element of a C*-algebra $A$
possesses a right and a left carrier projection in its bidual Banach space
$A^{**}$, a von Neumann algebra, and all the structural elements on Hilbert
$A$-modules can be canonically extended to the setting of Hilbert
$A^{**}$-modules, see the appendix and \cite{Pa1,Fr1} for details.

A subset $\{ x_j : j \in {\mathbb J} \}$ of $\mathcal H$  is said to be a
{\it generalized generating set} of the Hilbert $A$-module $\{ {\mathcal H},
\langle .,. \rangle \}$ if the $A$-linear hull of $\{ x_j : j \in {\mathbb J}
\}$ (i.e.~the set of all finite $A$-linear combinations of elements of this
set) is dense with respect to the topology induced by the semi-norms
$\{ |f(\langle .,. \rangle)|^{1/2} : f \in A^* \}$ in norm-bounded subsets of
$\mathcal H$. A generalized generating set is a {\it generalized Hilbert basis}
if its elements fulfil the conditions (i) and (ii) of the Hilbert basis
definition. The choice of the topology is motivated by its role in the
characterization of {\it self-dual} Hilbert C*-modules (i.e.~Hilbert
C*-modules $\mathcal H$ for which the Banach $A$-module ${\mathcal H}'$ of
all bounded $A$-linear maps $r: {\mathcal H} \to A$ coincides with
$\mathcal H$, \cite[Th.~6.4]{Frank:93}) and by the role of the weak*
topology for the characterization of Hilbert W*-modules and their special
properties (cf.~\cite{Pa1,Fr1} and the appendix).
In general, we have to be very cautious with the use of a C*-theoretical
analogue of the concept of linear independence for C*-modules since subsets of
C*-algebras $A$ may contain zero-divisors.

\smallskip
We are especially interested in finitely and countably generated Hilbert
C*-modules over unital C*-algebras $A$. A Hilbert $A$-module $\{ {\mathcal H},
\langle .,. \rangle \}$ is {\it (algebraically) finitely generated} if there
exists a finite set $\{x_1, ... , x_n \}$ of elements of $\mathcal H$ such that
every element $x \in \mathcal H$ can be expressed as an $A$-linear combination
$x = \sum_{j=1}^n a_j x_j$ ($a_j \in A$). Note, that {\it topologically
finitely generated Hilbert C*-modules} form a larger class than algebraically
finitely generated Hilbert C*-modules, cf.~Example \ref{ex214}. We classify the
non-algebraic topological case as belonging to the countably generated case
that is described below.

Algebraically finitely generated Hilbert $A$-modules over unital C*-algebras
$A$ are precisely the finitely generated projective $A$-modules in a pure
algebraic sense, cf.~\cite[Cor.~15.4.8]{Wegge-Olsen}.
Therefore, any finitely generated Hilbert $A$-module can be represented as an
orthogonal summand of some finitely generated free $A$-module $A^N = A_{(1)}
\oplus ... \oplus A_{(N)}$ consisting of all $N$-tuples with entries from $A$,
equipped with coordinate-wise operations and the $A$-valued inner product
$\langle (a_1, ..., a_N),(b_1, ... ,b_N) \rangle = \sum_{j=1}^N a_jb_j^*$.
The finitely generated free $A$-modules $A^N$ can be alternatively represented
as the algebraic tensor product of the C*-algebra $A$ by the Hilbert space
${\mathbb C}^N$.

Finitely generated Hilbert C*-modules have analogous properties to Hilbert
spaces in many ways. For example, they are self-dual, any bounded C*-linear
operator between two of them has an adjoint operator, and if they appear as
a Banach $A$-submodule of another Hilbert $A$-module we can always separate
them as an orthogonal summand therein.

\smallskip
The second and more delicate class of interest are the countably
generated Hilbert C*-modules over unital C*-algebras $A$. A Hilbert $A$-module
is {\it countably generated} if there exists a countable set of generators.
By G.~G.~Kasparov's Stabilization Theorem  \cite[Th.~1]{Kas} any countably
generated Hilbert $A$-module $\{ {\mathcal H}, \langle .,. \rangle \}$ over a
($\sigma$-)unital C*-algebra $A$ can be represented as an orthogonal summand
of the standard Hilbert $A$-module $l_2(A)$ defined by
\begin{equation} \label{def-l2}
l_2(A) = \left\{ \{ a_j : j \in {\mathbb N} \} \, : \, \sum_j a_ja_j^* \,\,\,
{\rm converges} \,\, {\rm in} \,\,\, \|.\|_A \right\} \, ,
\,
\langle \{ a_j \}, \{ b_j \} \rangle = \sum_j a_j b_j^* \, ,
\end{equation}
in such a way that its orthogonal complement is isomorphic to $l_2(A)$
again (in short: $l_2(A) \cong {\mathcal H} \oplus l_2(A)$). Often there exist
also different more complicated embeddings of $\mathcal H$ into $l_2(A)$.

As a matter of fact countably generated Hilbert C*-modules possess still the
great advantage that they are unitarily isomorphic as Hilbert $A$-modules iff
they are isometrically isomorphic as Banach $A$-modules, iff they are simply
bicontinuously isomorphic as Banach $A$-modules, \cite[Th.~4.1]{Frank:93}. So
we can omit the indication what kind of $A$-valued inner product on $\mathcal H$
will be considered because any two $A$-valued inner products on $\mathcal H$
inducing equivalent norms to the given one are automatically unitarily
isomorphic.

Countably generated Hilbert $A$-modules $\mathcal H$ are self-dual in only a
few cases. A large class consists of (countably generated) Hilbert $A$-modules
over finite-dimensional C*-algebras $A$ (i.e.~matrix algebras). However,
$l_2(A)$ is self-dual if and only if $A$ is finite-dimensional (\cite{Fr1}),
so further examples depend strongly on the special structure of the module
under consideration. In general, the $A$-dual Banach $A$-module $l_2(A)'$ of
$l_2(A)$ can be identified with the set
\[
l_2(A)'= \left\{ \{ a_j : j \in {\mathbb N} \} \, : \, \sup_{N \in {\mathbb N}}
\left\|\sum_{j=1}^N a_ja_j^* \right\|_A < \infty  \right\} \, .
\]
Every Hilbert C*-module possesses a standard isometric embedding into its
C*-dual Banach $A$-module via the $A$-valued inner product $\langle .,. \rangle$
defined on it varying the second argument of $\langle .,. \rangle$ over all
module elements.
The $A$-valued inner product on $l_2(A)$ can be continued to an $A$-valued
inner product on $l_2(A)'$ iff $A$ is a monotone sequentially complete
C*-algebra (e.g.~W*-algebra, monotone complete C*-algebra and little beyond).
So, for general considerations we have to face that ${\mathcal H} \not\equiv
{\mathcal H}'$ is the standard situation.

As a consequence of the lack of a general analogue of Riesz' theorem for
bounded module $A$-functionals on countably generated Hilbert $A$-modules
non-adjointable operators on $l_2(A)$ may exist, and they exist in fact
for every unital, infinite-dimensional C*-algebra $A$, cf.~\cite[Th.~4.3]{Fr1}
\cite[Cor.~5.6, Th.~6.6]{Frank:93}.
Furthermore, Banach C*-submodules can be either orthogonal summands,
or direct summands in a topological way only, or even they can lack the direct
summand property in any sense, cf.~\cite[Prop.~5.3]{Frank:93}. There are some
further surprising situations in Hilbert C*-module theory which cannot happen
in Hilbert space theory. Due to their minor importance for our considerations
we refer the interested reader to the standard reference sources on Hilbert
C*-modules \cite{Pa1,Rie,Kas,Jensen/Thomsen,Lance2,Wegge-Olsen,RaeWil,Bl1,Fr7}.

%\smallskip
%The existence of a finite or countable set of generators for a Hilbert
%C*-module $\{ {\mathcal H}, \langle .,. \rangle \}$ can be alternatively
%characterized in terms of the corresponding C*-algebra of "compact" operators
%${\rm K}_A({\mathcal H})$ on $\mathcal H$. This C*-algebra is defined as the
%operator norm closure of the linear span of the set
%$\{ \theta_{x,y} \, : \, \theta_{x,y}(z)= \langle z,x \rangle y \,\,\, {\rm
%for} \,\,\, z \in {\mathcal H} \}$ of elementary operators on $\mathcal H$.
%Note in passing that the property of an operator to be ''compact'' strongly
%depends on the choice of the $A$-valued inner product on $\mathcal H$, however
%for countably generated Hilbert C*-modules the C*-algebras of ''compact''
%operators arising from different but equivalent inner products are
%automatically $*$-isomorphic, cf.~\cite[Th.~4.1]{Frank:93}. The C*-algebra of
%''compact'' operators ${\rm K}_A({\mathcal H})$ on a Hilbert $A$-module
%$\mathcal H$ contains an identity if and only if the $A$-module $\mathcal H$
%is algebraically finitely generated, \cite[Exerc.~15.O]{Wegge-Olsen}. It is
%$\sigma$-unital if and only if the Hilbert $A$-module $\mathcal H$ is countably
%generated, \cite[Cor.~1.1.25]{Jensen/Thomsen}. Both criteria hold without
%reference to a concrete pre-selected $A$-valued inner product on $\mathcal H$
%within the unitary equivalence class.

\smallskip
If we consider finitely generated Hilbert C*-modules we do in general not
have any concept of a dimension since generating sets of elements can be
generating and irreducible at the same time and may, nevertheless, contain
different numbers of elements.

\begin{example}  \label{ex5} {\rm
   Let $A$ be the W*-algebra of all bounded linear operators on the separable
   Hilbert space $l_2$. Since the direct orthogonal sum of two copies of $l_2$
   is unitarily isomorphic to $l_2$ itself the projections $p_1$, $p_2$ to
   them are similar to the identity operator. Denote by $u_1$, $u_2$ the
   isometries realizing this similarity, i.e. $u_iu_i^*=1_A$, $u_i^*u_i=p_i$
   for $i=1,2$.  We claim that the Hilbert $A$-modules ${\mathcal H}_1=A$ and
   ${\mathcal H}_2=A^2$ are canonically isomorphic.
   \newline
   Indeed, the mapping $T : A \to A^2$, $T(a)=(au_1^*,au_2^*)$ (where
   $T^{-1}(c,d)=cu_1+du_2$) with $a,c,d \in A$ realize this unitary isomorphism.
   Consequently, ${\mathcal H}_1$ possesses two $A$-linearly independent sets of
   generators $\{ 1_A \}$ and $\{ u_1,u_2 \}$ with a different number of
   elements. Moreover, the ''magic'' formula (\cite[Cor.~1.2, (iii)]{HanLarson})
   $\sum \langle x_j,x_j \rangle = {\rm dim}(H)$ for frames $\{ x_j \}$ in
   Hilbert spaces $H$ does not work any longer: $1_A \cdot 1_A^* = 1_A$ and
   $u_1u_1^* + u_2u_2^* = 2 \cdot 1_A$.
   \newline
   In fact, for this C*-algebra $A$ the Hilbert $A$-module $A$ is unitarily
   isomorphic to $A^N$ for every $N \in {\mathbb N}$, $N \geq 0$, and the sum
   realizes the values $N \cdot 1_A$ for appropriate bases consisting of
   partial isometries.   }
\end{example}

What seems to be bad from the point of view of dimension theory of Hilbert
spaces sounds good from the point of view of frames. Normalized tight frames
of finitely gene\-rated Hilbert spaces have a number of elements that is
greater-equal the dimension of the Hilbert space under consideration,
cf.~\cite[Example ${\rm A}_1$]{HanLarson}. The number of elements of a frame
has never been an invariant of the Hilbert space. Therefore, the phenomena
fits into the already known picture quite well. What is more, concepts like
equivalence or similarity always compare frames with the same number of
elements, i.e.~are already restrictive in Hilbert space theory.

\smallskip
Concluding our introductory remarks about Hilbert C*-modules we want to fix
two further denotations. The set of all bounded $A$-linear operators on
$\mathcal H$ is denoted by ${\rm End}_A({\mathcal H})$, whereas the subset of
all adjointable bounded $A$-linear operators is denoted by
${\rm End}_A^*({\mathcal H})$.

%%%%%%%%%%%%%%%%%%%%%%%%%%%%%%%%%%%%%%%%%%%%%%%%%%%%%%%%%%%%%%%%%%%%%%%%%%%%%%
\section{Basic definitions}

The theory presented in this section is built up from basic principles of
functional analysis. We adopt the geometric dilation point of view of Deguang
Han and David R.~Larson in \cite{HanLarson}. To circumvent uncountable sets we
restrict ourself to countable frames. Uncountable frames cannot appear in
finite-dimensional Hilbert spaces (see Proposition \ref{prop-magic}) or
in separable Hilbert spaces (because of spectral theory), however they may
arise for e.g.~Hilbert ${\rm C}(X)$-modules since the underlying compact
Hausdorff space $X$ can be very complicated.

\begin{definition}  {\rm
    Let $A$ be a unital C*-algebra and ${\mathbb J}$ be a finite or countable
    index subset of $\mathbb N$.
    A sequence $\{ x_j : j \in {\mathbb J} \}$ of elements in a Hilbert
    $A$-module $\mathcal H$ is said to be a {\it frame} if there are real
    constants $C,D > 0$ such that
    \begin{equation} \label{ineq-frame}
      C \cdot \langle x,x \rangle \leq
      \sum_{j=1}^\infty \langle x,x_j \rangle \langle x_j,x \rangle \leq
      D \cdot \langle x,x \rangle
    \end{equation}
    for every $x \in \mathcal H$. The optimal constants (i.e.~maximal for $C$ and
    minimal for $D$) are called {\it frame bounds}. The frame $\{ x_j :
    j \in {\mathbb J} \}$ is said to be a {\it tight frame} if $C=D$, and said
    to be {\it normalized} if $C=D=1$. We consider {\it standard} (normalized
    tight) frames in the main for which the sum in the middle of the inequality
    (\ref{ineq-frame}) always converges in norm.

    A sequence $\{ x_j : j \in {\mathbb J}  \}$ is said to be a {\it
    (generalized) Riesz basis} if $\{ x_j : j \in {\mathbb J}  \}$ is a frame
    and a generalized generating set with one additional property: $A$-linear
    combinations $\sum_{j \in S} a_jx_j $ with coefficients $\{ a_j : j \in S \}
    \in A$ and $S \in {\mathbb J}$ are equal to zero if and only if in
    particular every summand $a_jx_j$ equals zero, $j \in S$. We call a
    sequence $\{ x_j : j \in {\mathbb J} \}$ in $\mathcal H$ a {\it standard
    Riesz basis for $\mathcal H$} if $\{ x_j : j \in {\mathbb J} \}$ is a frame
    and a generating set with the mentioned above uniqueness property for the
    representation of the zero element.
    An {\it inner summand of a standard Riesz basis} of a Hilbert $A$-module
    $\mathcal L$ is a sequence $\{ x_j : j \in {\mathbb J}  \}$ in a Hilbert
    $A$-module $\mathcal H$  for which there exists a second sequence $\{ y_j :
    j \in {\mathbb J}  \}$ in another Hilbert $A$-module $\mathcal K$ such
    that ${\mathcal L} \cong {\mathcal H} \oplus {\mathcal K}$ and the sequence
    consisting of the pairwise orthogonal sums $\{ x_j \oplus y_j : j \in
    {\mathbb J} \}$ in the Hilbert $A$-module ${\mathcal H} \oplus {\mathcal K}$
    is the original standard Riesz basis of $\mathcal L$.    }
\end{definition}

Since the set of all positive elements of a C*-algebra has the structure
of a cone the property of a sequence being a frame does not depend on
the sequential order of its elements. Consequently, we can replace the
ordered index set ${\mathbb J} \subseteq {\mathbb N}$ by any countable index
set ${\mathbb J}$ without loss of generality. We do this for further purposes.

In Hilbert space theory a Riesz basis is sometimes defined to be a basis
arising as the image of an orthonormal basis by an invertible linear operator.
Since the concept of orthonormality cannot be transfered one-to-one to the
theory of Hilbert C*-modules the suitable generalization of this statement
needs to clarify this. Especially, the more complicated inner structure
of C*-algebras $A$ in comparison to the field of complex numbers ${\mathbb C}$
has to be taken into account. We will formulate an analogous result as
Corollary \ref{cor-Riesz} below. The other way around standard Riesz bases can
be characterized as frames $\{ x_i : i \in \mathbb J \}$ such that the
$A$-module generated by one single element $x_j$ of the frame has always only
a trivial intersection with the norm-closed $A$-linear span of the other
elements $\{ x_i : i \not= j \}$.

\smallskip
The definition above has some simple consequences. A set $\{ x_j : j \in
{\mathbb J} \}$ is a normalized tight frame if and only if
the equality
  \begin{equation} \label{eq-ntframe}
    \langle x,x \rangle =
       \sum_{j \in {\mathbb J}} \langle x,x_j \rangle \langle x_j,x \rangle
  \end{equation}
holds for every $x \in \mathcal H$. Note that this sum can fail to converge
uniformly in $A$, however the sum always converges in $A$ with respect to the
weak topology induced by the dual space $A^*$ of $A$ (cf.~Example \ref{ex1}
below).

Furthermore, the norms of the elements of a frame are always uniformly bounded
by the square root of the upper frame bound $D$. To see that consider the
chain of inequalities
  \[
    \langle x_k,x_k \rangle^2
    \leq \sum_{j \in {\mathbb J}} \langle x_k,x_j \rangle \langle x_j,x_k \rangle
    \leq D \cdot \langle x_k,x_k \rangle
  \]
that is valid for every $k \in {\mathbb J}$. Taking the norms on both sides
the inequality is preserved.

\begin{proposition}  \label{on-basis2}
   Let $A$ be a C*-algebra and $\mathcal H$ be a finitely or countably generated Hilbert
   $A$-module.
   \begin{itemize}
   \item[(i)] If an orthogonal Hilbert basis $\{ x_j : j \in {\mathbb J} \}$
      of $\mathcal H$ is a standard normalized tight frame then the values
      $\{ \langle x_j,x_j \rangle : j \in {\mathbb J} \}$ are all non-zero
      projections.
   \item[(ii)] Conversely, every standard normalized tight frame $\{ x_j :
      J \in {\mathbb J} \}$ of $\mathcal H$ for which the values $\{ \langle
      x_j,x_j \rangle : j \in {\mathbb J} \}$ are non-zero projections is an
      orthogonal Hilbert basis of $\mathcal H$.
   \end{itemize}
   In general, the inequality $\langle x_j,x_j \rangle \leq 1_A$ holds for every
   element $x_j$ of normalized tight frames $\{ x_j : J \in {\mathbb J} \}$ of
   $\mathcal H$.
\end{proposition}

\begin{proof}
Fix an orthogonal Hilbert basis $\{ x_j : j \in {\mathbb J} \}$ of $\mathcal H$.
Consider norm-convergent sums $x=\sum_j a_j x_j \in \mathcal H$ for suitably
selected sequences $\{ a_j : j \in {\mathbb J} \} \in A$.
If the Hilbert basis of $\mathcal H$ is a normalized tight frame then the
equality
   \begin{eqnarray*} 
      \sum_{j \in {\mathbb J}} a_j \langle x_j,x_j \rangle a_j^*  & = &
               \left\langle \sum_{j \in {\mathbb J}} a_j x_j ,
               \sum_{k \in {\mathbb J}} a_k x_k  \right\rangle 
           =   \langle x,x \rangle \\
         & = & \sum_{j \in {\mathbb J}} \langle x,x_j \rangle \langle x_j,x
               \rangle 
           =   \sum_{j \in {\mathbb J}} \left\langle \sum_{k \in {\mathbb J}}
               a_k x_k , x_j \right\rangle \left\langle x_j ,
               \sum_{l \in {\mathbb J}} a_l x_l  \right\rangle \\
         & = & \sum_{j \in {\mathbb J}} \langle a_jx_j,x_j \rangle \langle x_j,
               a_jx_j \rangle 
           =   \sum_{j \in {\mathbb J}} a_j \langle x_j,x_j \rangle^2 a_j^*
   \end{eqnarray*}
is valid for every admissible choice of the coefficients $\{ a_j : j \in
{\mathbb J} \} \in A$. In particular, one admissible selection is $a_i=1_A$
and $a_j=0_A$ for each $j \not= i$, $i \in {\mathbb J}$ fixed. For this setting
we obtain $0 \not= \langle x_i,x_i \rangle = \langle x_i,x_i \rangle^2$ since
$x_i \not= 0$ by supposition.

The converse conclusion is also a simple calculation. If $\{ x_j : j \in
{\mathbb J} \}$ is a standard normalized tight frame, then (\ref{eq-ntframe})
implies
   \begin{equation*}  
      0 \leq \sum_{j \not= i} \langle x_i,x_j \rangle\langle x_j,x_i \rangle
          =           \langle x_i,x_i \rangle - \langle x_i,x_i \rangle^2 \, .
   \end{equation*}
Therefore, $\langle x_j,x_j \rangle \leq 1_A$ for every $j \in {\mathbb J}$ by
spectral theory. Now, if some element $x_i \not= 0$ happens to admit a
projection as the inner product value $\langle x_i,x_i \rangle$, then
$0 = \sum_{j \not= i} \langle x_j,x_i \rangle\langle x_i,x_j \rangle$,
i.e. $\langle x_j,x_i \rangle$ for any $j \not= i$ by the positivity of
the summands. In other words, the element $x_i$ must be orthogonal to all
other elements $x_j$, $j \not= i$, of that normalized tight frame. Consider
a decomposition of the zero element in the special form $0=\sum_j a_jx_j$ for
suitably selected  coefficients $\{ a_j : j \in {\mathbb J} \} \subset A$.
Since
\[
   0  =  \left\langle \sum_{j \in {\mathbb J}} a_jx_j, \sum_{k \in {\mathbb
            J}} a_kx_k \right\rangle 
      =  \sum_{j \in {\mathbb J}} \langle a_jx_j,a_jx_j \rangle
\]
and since the sum at the right end is a sum of positive summands we arrive at
$a_jx_j = 0$ for every $j \in \mathbb J$. Thus, a standard normalized tight
frame $\{ x_j : j \in {\mathbb J} \}$ for which the values $\{ \langle x_j,x_j
\rangle : j \in {\mathbb J} \}$ are non-zero projections is an orthogonal
Hilbert basis of $\mathcal H$.
\end{proof}

As in the Hilbert space situation we would like to establish that standard
Riesz bases that are normalized tight frames have to be orthogonal Hilbert
bases with projections as the values of the inner products with equal
basis element entries. This requires some more work than expected and has to
be postponed until we derive the reconstruction formula, cf.~Corollary
\ref{on-basis}.

\begin{lemma}
   Let $A$ be a unital C*-algebra. For some element $x$ of a Hilbert C*-module
   $\{ {\mathcal H}, \langle .,. \rangle \}$ the elementary ''compact''
   operator $\theta_{x,x}$ mapping $y \in \mathcal H$ to $\langle y,x \rangle x$
   is a projection if and only if
   $x=\langle x,x \rangle x$, if and only if $\langle x,x \rangle$ is a
   projection. In this case the elements of $Ax \subseteq \mathcal H$ can be
   identified with the elements of the ideal $A\langle x,x \rangle \subseteq A$.
   If for two orthogonal elements $x,y \in \mathcal H$ with $x = \langle x,x
   \rangle x$, $y = \langle y,y \rangle y$ the equality $\langle x,x \rangle =
   \langle y,y \rangle$ holds additionally, then the projections $\theta_{x,x}$
   and $\theta_{y,y}$ are similar in the sense of Murray-von Neumann, where
   the connecting partial isometry is $\theta_{x,y}$.
\end{lemma}

The statement can be verified by elementary calculations and, thus, a proof
is omitted.

Since there exist unital C*-algebras $A$ such that the monoid of all finitely
generated projective $A$-modules with respect to orthogonal sums does
not possess the cancelation property, in some situations orthogonal Hilbert
or Riesz bases may not exist. Examples can be found in sources about operator
$K$-theory of C*-algebras, cf.~\cite{Wegge-Olsen}.
Also, for unital C*-algebras $A$ with an extremly small subset of orthogonal
projections we are faced with countably generated Hilbert $A$-modules
$\mathcal H$ without any orthogonal Riesz basis. To give an example let
$A={\rm C}([0,1])$ be the C*-algebra of continuous function on the unit
interval and ${\mathcal H}={\rm C}_0((0,1])$ be the Hilbert
${\rm C}([0,1])$-module of all continuous functions on $[0,1]$ vanishing at
zero, equipped with the standard $A$-valued inner product. By the
Stone-Weierstrass theorem the set of functions $\{ t,t^2,...,t^n,... \}
\subset {\rm C}_0((0,1])$, $(t \in [0,1])$, possesses a ${\mathbb C}$-linear
hull that is norm-dense in ${\rm C}_0((0,1])$ and hence, $\mathcal H$ is a
countably generated Hilbert ${\rm C}([0,1])$-module. However, $\langle x,x
\rangle = \langle x,x \rangle^2$ for some $x \in {\rm C}_0((0,1])$ if and only
if $x=0$ since the only non-trivial projection $1_A \in A$ cannot be admitted
for inner product values of elements from $\mathcal H$. Nevertheless,
${\mathcal H}={\rm C}_0((0,1])$ has standard normalized tight frames as a
Hilbert ${\rm C}([0,1])$-module, see Example \ref{ex2} below.Also there exists
a trivial orthogonal Hilbert basis consisting of the single element $\{ t \}$.

\begin{example} \label{ex214} {\rm
  If $A$ is a unital C*-algebra and $\mathcal H$ is a countably generated
  Hilbert $A$-module then there may exist orthogonal Hilbert bases $\{ x_j \}$
  of $\mathcal H$ without the property $\langle x_i,x_i \rangle = \langle x_i,
  x_i \rangle^2$ for $j \in {\mathbb N}$. By Proposition \ref{on-basis2},
  these Hilbert bases are not frames.
  The roots of the problem behind this phenomenon lie in the difference between
  algebraically and topologically finite generatedness of Hilbert C*-modules.

  For example, set $A={\rm C}([0,1])$ to be the C*-algebra of all continuous
  functions on the unit interval and consider the set and Hilbert $A$-module
  ${\mathcal H}=l_2({\rm C}_0((0,1]))$ (cf.~(\ref{def-l2}) for the definition),
  where ${\rm C}_0((0,1])$ denotes the C*-subalgebra of all functions on
  $[0,1]$ vanishing at zero.
  The function $f(t)=t$ for $t \in [0,1]$ is topologically a single generator
  of ${\rm C}_0((0,1])$ by the Stone-Weierstra{\ss} theorem. The Hilbert
  $A$-module $\mathcal H$ is generated by the set $\{ f_i = (0_A, ..., 0_A,
  f_{(i)}, 0_A, ... ) : i \in {\mathbb N} \}$ of pairwise orthogonal elements
  of norm one.
  However, the inner product values of all these elements equal $f^2$ which is
  not a projection and the spectrum of which is not deleted away from zero.
  Therefore, the lower frame bound has to be zero.

  Looking for another orthogonal standard Riesz basis $\{ f_j : j \in \mathbb J
  \}$ of $\mathcal H$ we can only consider bases with two or more elements.
  However, $f_i \perp f_j$ always means that there exists a point $t_0 \in
  (0,1]$ such that $f_i \equiv 0$ for small $t \leq t_0$ and $f_j \equiv 0$ for
  small $t \geq t_0$. Taking into account orthogonality of these elements
  $\{ f_j \}$ every function in the norm-closed $A$-linear hull of them has to
  be zero at $t_0$ contradicting the assumptions. The only possible conclusion
  is the non-existence of any orthogonal standard Riesz basis of $\mathcal H$.
  We will see at Corollary \ref{cor-Riesz} that the existence of a standard
  Riesz basis of $\mathcal H$ would imply the existence of an orthogonal
  Hilbert basis for it that is a (standard) normalized tight frame at the same
  time. Therefore, $\mathcal H$ does even not possess any standard Riesz basis.
  }
\end{example}

In this place  we can state the following about standard Riesz bases of Hilbert
C*-modules (cf.~Corollary \ref{cor-Riesz}):

\begin{proposition}   \label{prop-orthoRiesz}
  Let $A$ be a unital C*-algebra and $\mathcal H$ be a countably or finitely
  generated Hilbert $A$-module. If $\mathcal H$ possesses an orthogonal
  standard Riesz basis then $\mathcal H$ possesses an orthogonal standard
  Riesz basis $\{ x_j : j \in {\mathbb J} \}$ with the property
  $\langle x_j,x_j \rangle = \langle x_j,x_j \rangle^2$ for any $j \in
  \mathbb J$, i.e.~an orthogonal Hilbert basis that is a standard
  normalized tight frame.
\end{proposition}

\begin{proof}
Suppose, $\mathcal H$ possesses an orthogonal standard Riesz basis
$\{ x_j \}$. That means, there are two constants $0 < C,D$ such that the
inequality $C \cdot \langle x_j,x_j \rangle \leq \langle x_j,x_j \rangle^2
\leq D \cdot \langle x_j,x_j \rangle$ is fulfilled for every $j \in {\mathbb J}$.
Obviously, $D=1$ since $\{ x_j \}$ is supposed to be a Hilbert basis and,
therefore, $\| x_j \| =1$ by one of the properties of Hilbert bases.
Considering the lower estimate with the constant $C$ spectral theory forces
the spectra of the elements $\{ \langle x_j,x_j \rangle \}$ to be uniformly
bounded away from zero by this constant $C$. Consequently, there are continuous
positive functions $\{f_j\}$ on the spectra of the elements $\{ \langle
x_j,x_j \rangle \}$ such that $f_j \langle x_j,x_j \rangle = (f_j \langle
x_j,x_j \rangle)^2$ and the restriction of these functions to the bounded away
from zero part of the spectra of $\{ \langle x_j,x_j \rangle \}$ equals one.
The new frame $\{ f_i^{1/2} x_j \}$ is normalized tight and orthogonal.
Moreover, it is standard since the spectra of the inner product values were
uniformly bounded away from zero.
\end{proof}

On the other hand, a frame may contain the zero element arbitrary often.
Indeed, if a normalized tight frame $\{ x_j : j \in {\mathbb J} \}$ has a
subsequence $\{ x_j : j \in {\mathbb I} \}$ that is a normalized tight frame for
$\mathcal H$, too, then (\ref{eq-ntframe}) implies
  \[
    \langle x_k, x_k \rangle =
    \sum_{j \in {\mathbb J}} \langle x_k,x_j \rangle \langle x_j,x_k \rangle =
    \sum_{i \in {\mathbb I}} \langle x_k,x_i \rangle \langle x_i,x_k \rangle  \, ,
  \]
i.e.~$\sum_{j \in {\mathbb J} \setminus {\mathbb I}} \langle x_k,x_j \rangle
\langle x_j,x_k \rangle = 0$. In case $k \in {\mathbb J} \setminus {\mathbb I}$
we obtain $\langle x_k, x_k \rangle^2=0$ and hence, $x_k=0$. Consequently,
normalized tight frames are maximal generating sets in some sense.

However, frames $\{ x_j : j \in \mathbb J \}$ may fail to meet the most
important property of a Hilbert basis of $\mathcal H$, nevertheless. As known
by examples of frames of two-dimensional Hilbert spaces $\mathcal H$ they may
contain to much elements to be a Hilbert basis of $\mathcal H$ since the
uniqueness of decomposition of elements $x \in {\mathcal H}$ as $x=\sum_j a_j
x_j$ for elements $\{ a_j : j \in {\mathbb J} \} \subset A$ may not be
guaranteed any longer (\cite[Example ${\rm A}_1$]{HanLarson}), in particular
the representation of the zero element can be realized as a sum of on-zero
summands.

\begin{definition}  {\rm
   Frames $\{ x_j : j \in {\mathbb J} \}$ and $\{ y_j : j \in {\mathbb J} \}$ of
   Hilbert $A$-modules $\mathcal H$ and $\mathcal K$, respectively, are
   {\it unitarily equivalent} if there is an $A$-linear unitary operator
   $U: {\mathcal H} \to {\mathcal K}$ such that $U(x_j)=y_j$ for every $j \in {\mathbb J}$.
   They are {\it similar} ({\it or isomorphic}) if the operator $U$ is merely
   bounded, adjointable, $A$-linear and invertible.  }
\end{definition}

We want to note that isomorphisms of frames are in general not invariant under
permutations, especially, if the frames contain the zero element. Moreover,
frames of different size in finitely generated Hilbert C*-modules cannot be
related by these concepts. To achieve sufficiently strong statements we
will not go into further modifications of similarity and isomorphism
concepts for frames.

%%%%%%%%%%%%%%%%%%%%%%%%%%%%%%%
\section{Examples of frames}

\begin{example}  \label{ex0}   {\rm
   Every sequence $\{ x_j : j \in \mathbb J \}$ of a finitely or countably
   generated Hilbert $A$-module for which every element $x \in \mathcal H$
   can be represented as $x=\sum_j \langle x,x_j \rangle x_j$ (in a probably
   weaker sense of series convergence than norm-convergence) is a normalized
   tight frame in $\mathcal H$. The decomposition of elements of $\mathcal H$
   is norm-convergent if and only if $\{ x_j : j \in \mathbb J \}$ is a
   standard normalized tight frame. Indeed,
   \begin{eqnarray*}
   \langle x,x \rangle & = &
      w-\lim_{n \to \infty}
      \left\langle \sum_{k=1}^n \langle x,x_k \rangle x_k, x \right\rangle \\
      & = &  w-\lim_{n \to \infty}
      \sum_{k=1}^n \langle \langle x,x_k \rangle x_k,x \rangle \\
      & = &  w-\lim_{n \to \infty}
      \sum_{k=1}^n \langle x,x_k \rangle \langle x_k,x \rangle  \, .
   \end{eqnarray*}
   }
\end{example}

\begin{example}  \label{ex-index}   {\rm
  Let $B$ be a unital C*-algebra and $E: B \to A \subseteq B$ be a conditional
  expectation on $B$. By Y.~Watatani $E$ is said to be algebraically of finite
  index if there exists a finite family $\{ (u_1,v_1), ..., (u_n,v_n) \}
  \subseteq B \times B$ that is called a {\it quasi-basis} such that
    \[
       x = \sum_i u_iE(v_ix) = \sum_i E(xu_i)v_i
    \]
  for every $x \in B$, cf. \cite[Def.~1.2.2]{Wata}. These expressions can be
  translated as decompositions of $B$ as a right/left finitely generated
  projective $A$-module, and it can be seen to be derived from an $A$-valued
  inner product on $B$ setting $\langle .,. \rangle =E(\langle .,. \rangle_B)$.
  We will see in section \ref{sec-dual} that the sets $\{ u_1, ... , u_n \}$
  and $\{ v_1, ... , v_n \}$ are dual to each other frames of $B$ as a
  finitely generated Hilbert $A$-module. Moreover, the setting $v_i=u_i^*$
  is the choice for the canonical dual of a normalized tight frame
  $\{ u_1, ... , u_n \}$, and such a choice can be made in every situation
  (see \cite[Lemma 2.1.6]{Wata}). The concept survives an extension to
  faithful bounded $A$-bimodule maps on $B$, \cite[Def.~1.11.2]{Wata}.

  To give a concrete example consider the matrix C*-algebra $B=M_n({\mathbb C})$
  and the norma\-lized trace $E={\rm tr}$ on it. The quasi-basis may be derived,
  for example, from a sequence of $n$ pairwise orthogonal minimal projections
  $\{ P_i : i = 1,...,n \} \in B$ and of a set of minimal partial isometries
  $\{ U_j : j=1, ... , n(n-1) \} \in B$ connecting them pairwise.
  As a special choice
  we could take the matrix-units, i.e.~all matrices with exactly one non-zero
  entry that equals one. Taking the selected $n^2$ elements of $B$ as the
  first part $\{ u_i : i=1,...,n \}$ of a suitable quasi-basis, and setting
  $v_i=u_i^*$, $i=1,...,n$, for the second part of it we obtain
    \[
      X = \sum_{i=1}^n u_i \cdot {\rm tr}(u_i^*X) =
          \sum_{i=1}^n {\rm tr}(Xu_i^*) \cdot u_i
    \]
  for $X \in B$ as desired. Of course, we have the special situation of $A =
  {\mathbb C}$, i.e.~a Hilbert space $\{ B , {\rm tr}(\langle .,. \rangle_B) \}$.

  Now, let $A$ be the subset of all diagonal matrices in $B$ and let $E$ be
  the mapping acting as the identity mapping on the diagonal of a matrix and
  as the zero mapping on off-diagonal elements. We can take the same
  quasi-basis as before, and we obtain
    \[
      X = \sum_{i=1}^n u_i \cdot E(u_i^*X) =
          \sum_{i=1}^n E(Xu_i^*) \cdot u_i
    \]
  for $X \in B$, again.

  Much more complicated examples are known for type II and III W*-factors,
  and for certain C*-algebras beyond the W*-class
  \cite{Jones,Kosaki,PP,KiFr98}  }
\end{example}

\begin{example}   \label{ex1}    {\rm
   Let $H$ be an infinite-dimensional Hilbert space and $\{ p_\alpha : \alpha
   \in I \}$ be a maximal set of pairwise orthogonal minimal orthogonal
   projections on $H$.
   Consider the C*-algebra $A={\rm B}(H)$ of all bounded linear operators on
   $H$ and the Hilbert $A$-modules ${\mathcal H}_1=A$ and ${\mathcal H}_2={\rm K}(H)$,
   where the latter consists of all compact operators on $H$. The set $\{
   p_\alpha \}$ is a normalized tight frame for both ${\mathcal H}_1$ and
   ${\mathcal H}_2$, however it is a non-standard one in the first case. Moreover,
   for this tight frame we obtain $\langle p_j, p_j \rangle = \langle p_j, p_j
   \rangle^2$ and $x = \sum_j \langle x,p_j \rangle p_j$ in the sense of
   w*-convergence in $A$. The frame is not a standard Riesz basis for
   ${\mathcal H}_1$ since it generates only ${\mathcal H}_2$ by convergence
   in norm. Note, that the frame can contain uncountably many elements.  }
\end{example}

The structural obstacle behind this phenomenon is order convergence.
Infinite-dimen\-sio\-nal C*-algebras $A$ can possess sequences of pairwise
orthogonal positive elements the sum of which converges in order inside $A$,
but not in norm. They may cause this kind of non-standard normalized tight
frames. Since the structure of the basic C*-algebra $A$ may be very
complicated containing monotone complete and non-complete blocks we have to
try to circumvent this kind of situation in our first attempt to generalize
the theory.
Otherwise, Theorem \ref{th-reconstr} can be only formulated for
self-dual Hilbert $A$-modules over monotone complete C*-algebras $A$ since
only for this class of Hilbert C*-modules the $A$-valued inner product can
be canonically continued to an $A$-valued inner product on the $A$-dual
Banach $A$-module of a given Hilbert $A$-module. The disadvantage consists
in the small number of examples covered by this setting, most of them far
from being typical. The other way out of the situation would be a switch
to general Banach $A$-module theory without any inner product structures.
This is surely possible but technically highly complicated. So we will
restrict ourself to standard frames for the time being.

\begin{example} \label{ex2}   {\rm
   Let $A$ be the C*-algebra of all continuous functions on the unit interval.
   Let $\mathcal H$ be the set of all continuous functions on $[0,1]$
   vanishing at zero. The set $\mathcal H$ is a countably generated Hilbert
   $A$-module by the Stone-Weierstrass theorem (take e.g.~the functions
   $\{t,t^2,t^3,... \}$ as a set of generators). The $A$-valued inner product
   on $\mathcal H$ is defined by the formula $\langle f,g \rangle = fg^*$.
   As already discussed this Hilbert $A$-module does not contain any orthogonal
   Riesz basis.

   However, $\mathcal H$ possesses standard normalized tight frames.
   The following set of elements of $\mathcal H$ forms one:
   \[
   x_j(t) = \left\{ \begin{array}{c@{\quad:\quad}l}
                  \sqrt{j(j+1)t-j} & t \in [(j+1)^{-1},j^{-1}] \\
                  \sqrt{-j(j-1)t+j} & t \in [j^{-1},(j-1)^{-1}] \\
                  0 & {\rm elsewhere}
                  \end{array} \right.  \,\,\qquad {\rm for} \quad j > 1 \, ,
   \]
   \[
   x_1(t) = \left\{ \begin{array}{c@{\quad:\quad}l}
                  \sqrt{2t-1} & t \in [1/2,1] \\
                  0 & t \in [0,1/2]
                  \end{array} \right.
   \]
   It is not a frame for the (singly generated) Hilbert $A$-module $A$ itself
   since the constant $C$ of inequality (\ref{ineq-frame}) has to be
   zero for this extended Hilbert $A$-module (look at $t=0$ for functions $f$
   with $f(0) \not= 0$). Adding a further element $x_0=f$ with $f(0) \not= 0$
   to the sequence under consideration we obtain a frame for the Hilbert
   $A$-module $A$, however not a tight one since $\, \max C = | f(0) |^2 \,$
   and $\, \min D = 1 + \max | f(t) |^2$.      }
\end{example}

%\begin{example}    {\rm
%   Let $A = l_\infty$. Denote the set of all complex-valued sequences
%   conver\-ging to zero by $c_0$.
%   Consider the Hilbert $A$-module ${\mathcal H}=l_\infty \oplus c_o$ and
%   its submodule ${\mathcal K} = \{ (i,i) : i \in c_0 \}$. The $A$-valued inner
%   product is defined in the standard manner. The sequence of ordered pairs
%   \[
%   x_1 = (1_A,0_A) \, , \, \, x_j = (0_A,(0,0, ... ,0,1_{(j-1)},0, ...))
%   \quad (j \in {\mathbb N})
%   \]
%   is a standard Riesz basis of $\mathcal H$. Since $\mathcal K$ is not an orthogonal
%   summand of $\mathcal H$, but a topological direct summand (with complement
%   $\{ (a,0_A) : a \in A \}$) the respectively projected to $\mathcal K$
%   sequence $\{ P(x_j) \} \in \mathcal K$ is not a normalized tight frame for
%   $\mathcal K$. However, this sequence is still a tight frame of $\mathcal K$
%   with $C=D=2$. }
%\end{example}

\begin{example}  \label{ex4} {\rm
  After these unusual examples we want to indicate
  good classes of frames for every finitely and countably generated Hilbert
  $A$-module $\mathcal H$ over a unital C*-algebra $A$. In fact, there is an
  abundance of standard normalized tight frames in each finitely or countably
  generated Hilbert $A$-module:
  recall that the standard Hilbert $A$-modules $A^N$ $(N \in {\mathbb N})$ and
  $l_2(A)$ have unitarily isomorphic representations as (normed linear space)
  tensor products of the C*-algebra $A$ and the Hilbert spaces ${\mathbb C}^N$
  $(N \in {\mathbb N})$ and $l_2({\mathbb C})$, respectively. Simply set the
  $A$-valued inner product to
    \[
    \langle a \otimes h, b \otimes g \rangle = ab^* \langle h,g \rangle_H
    \]
  for $a,b \in A$ and $g,h$ from the appropriate Hilbert space $H$. In fact,
  the algebraic tensor product $A \odot l_2$ needs completion with respect to
  the arising Hilbert norm to establish the unitary isomorphism.

  Using this construction every frame $\{ x_j \}$ of the Hilbert space $H$
  induces a corresponding frame $\{ 1_A \otimes x_j \}$ in $A^N$ $(N \in
  {\mathbb N})$ or $l_2(A)$. The properties to be tight or (standard) normalized
  tight transfer. Non-standard normalized tight frames in Hilbert C*-modules
  cannot arise this way.

  To find frames in arbitrary finitely or countably generated Hilbert
  C*-modules over unital C*-algebras $A$ recall that every such Hilbert
  $A$-module $\mathcal H$ is an orthogonal summand of $A^N$ $(N \in {\mathbb N})$
  or $l_2(A)$, respectively (see section one). So there exists an orthogonal
  projection $P$ of $A^N$ or $l_2(A)$ onto this embedding of $\mathcal H$. The
  next fact to show is that any orthogonal projection of an orthonormal Riesz
  basis of $A^N$ or $l_2(A)$ is a standard normalized frame of the range
  $\mathcal H$ of $P$.

  Denote the standard Riesz basis of $A^N$ or $l_2(A)$ by $\{ e_j \}$ and
  the elements of the resulting sequence $\{ P(e_j) \}$ by $x_j$, $j \in {\mathbb N}$.
  For every $x \in \mathcal H$ we have
  \[
  \langle x,x \rangle = \sum_j \langle x,e_j \rangle \langle e_j,x \rangle \quad
  , \quad
  x= \sum_j \langle x,e_j \rangle e_j    \, .
  \]
  Applying the projection $P$ to the decomposition of $x$ with respect to the
  orthonormal basis $\{ e_j \}$ we obtain $x = \sum_j \langle x,x_j \rangle x_j$
  since $x=P(x)$, $x_j=P(e_j)$ and $\langle x,e_j \rangle=\langle x,x_j \rangle$
  for $j \in {\mathbb N}$. By Example \ref{ex0} the sequence
  $\{ x_j \}$ becomes a standard normalized tight frame of $\mathcal H$.

  This formula $x = \sum_j \langle x,x_j \rangle x_j$ is called {\it the
  reconstruction formula of a frame} in Hilbert space theory.
  The remaining point is to show that {\it every} standard normalized tight
  frame of finitely and countably generated Hilbert $A$-modules over unital
  C*-algebras $A$ arises in this way, see Theorem \ref{th-reconstr} below
  (and even non-standard ones, see section eight).    }
\end{example}

%%%%%%%%%%%%%%%%%%%%%%%%%%%%%%%%%%%%%%%%%%%%%%%%%%%%%%%%

\section{Frame transform and reconstruction formula}

This section is devoted to the key result that allows all the further
developments we could work out. We found that for unital C*-algebras $A$
the frame transform operator related to a standard (normalized tight) frame
in a finitely or countably generated Hilbert $A$-module is adjointable in
every situation, and that the reconstruction formula holds.
Moreover, the image of the frame transform is an orthogonal summand of
$l_2(A)$. The proof is in crucial points different from that one for Hilbert
spaces since these properties of the frame transform are not guaranteed by
general operator and submodule theory. At the opposite, the results are rather
unexpected in their generality to hold and have to be established by
non-traditional arguments. For the Hilbert space situation we refer to
\cite[Prop.~1.1]{HanLarson} and \cite[Th.~2.1, 2.2]{Holub:94}.

\begin{theorem}  \label{th-reconstr} {\rm (frame transform and reconstruction
                                                                     formula)}
\newline
  Let $A$ be a unital C*-algebra, $\{ {\mathcal H}, \langle .,. \rangle \}$ be a
  finitely or countably generated Hilbert $A$-module.
  Suppose that $\{ x_n : n \in {\mathbb J} \}$ is a standard normalized tight
  frame for $\mathcal H$. Then the corresponding frame transform $\theta:
  {\mathcal H} \to l_2(A)$ defined by $\theta(x)=\{ \langle x,x_n \rangle
  \}_{n \in {\mathbb J}}$ for $x \in {\mathcal H}$ possesses an adjoint
  operator and realizes an isometric embedding of $\mathcal H$ onto an
  orthogonal summand of $l_2(A)$. The adjoint operator $\theta^*$ is surjective
  and fulfills $\theta^*(e_n)=x_n$ for every $n \in {\mathbb J}$.
  Moreover, the corresponding orthogonal projection $P: l_2(A) \to \theta
  ({\mathcal H})$ fulfills $P(e_n) \equiv \theta(x_n)$ for the standard orthonormal
  basis $\{ e_n = (0_A, ... ,0_A,1_{A,(n)},0_A, ... ) : n \in {\mathbb J} \}$ of
  $l_2(A)$. For every $ x \in {\mathcal H}$ the decomposition $x=\sum_i \langle
  x,x_i \rangle x_i$ is valid, where the sum converges in norm.
  \newline
  The frame $\{ x_n \}$ is a set of module generators of the Hilbert $A$-module
  $\mathcal H$. If the frame is not a Riesz basis then the frame elements do
  not form an $A$-linearly independent set of elements. The operator equality
  ${\rm id}_{{\mathcal H}} = \sum_i \theta_{x_i,x_i}$ is fulfilled in the sense
  of norm-convergence of the series $\sum_i \theta_{x_i,x_i}(x)$ to $x \in
  \mathcal H$.
\end{theorem}

\begin{proof}
Since the sequence $\{ x_j : j \in {\mathbb J} \}$ is a standard
normalized tight frame in $\mathcal H$ the frame operator is correctly defined
and the equality
  \[
    \langle \theta(x), \theta(x) \rangle_{l_2} =  \sum_{j \in {\mathbb J}}
    \langle x,x_j \rangle_{{\mathcal H}} \langle x_j,x \rangle_{{\mathcal H}} =
    \langle x,x \rangle_{{\mathcal H}}
  \]
holds for any $x \in \mathcal H$. Moreover, the image of $\theta$ is closed
because $\mathcal H$ is closed by assumption. So, $\theta$ is an isometric
$A$-linear embedding of $\mathcal H$ into $l_2(A)$ with norm-closed image.

\smallskip
To calculate the values of the adjoint operator $\theta^*$ of $\theta$
consider the equality
  \begin{eqnarray*}
    \langle \theta(x),e_i \rangle_{l_2(A)} & = &
    \left\langle \sum_k \langle x,x_k \rangle_{{\mathcal H}} e_k,
                                                e_i \right\rangle_{l_2(A)} \\
    & = & \sum_k \langle x,x_k \rangle_{{\mathcal H}} \langle e_k,
                                                e_i \rangle_{l_2(A)} \\
    & = & \langle x,x_i \rangle_{{\mathcal H}}
  \end{eqnarray*}
which is satisfied for every $x \in \mathcal H$, every $i \in \mathbb J$.
Consequently, $\theta^*$ is at least defined for the elements of the selected
orthonormal Riesz basis $\{ e_j : j \in \mathbb J \}$ of $l_2(A)$ and takes
the values $\theta^*(e_j)=x_j$ for every $j \in {\mathbb J}$.
Since the operator $\theta^*$ has to be $A$-linear by definition we can extend
this operator to the norm-dense subset of all finite $A$-linear combinations
of the elements of the selected basis of $l_2(A)$.

Furthermore, we are going to show that $\theta^*$ is bounded. To see this
consider the bounded $A$-linear mapping $\langle \theta(\cdot ), y \rangle$
from $l_2(A)$ to $A$ defined for any $y \in l_2(A)$. The inequality
  \begin{eqnarray*}
      \| \langle x,\theta^*(y) \rangle_{\mathcal H} \|_A
      & = &  \| \langle \theta(x),y \rangle_{l_2(A)} \|_A \\
      & \leq &  \| \theta \| \| x \| \| y \|
  \end{eqnarray*}
is valid for any $y$ that is an element of the domain of $\theta^*$ 
and for any $x \in \mathcal H$ by the general Cauchy-Schwarz inequality for
Hilbert C*-modules. Taking the supremum over the set $\{ x \in
{\mathcal H} : \| x \| \leq 1 \}$ of both the sides of the inequality we get
    \[
      \| \theta^*(y) \|_{{\mathcal H}'} =
      \| \langle \cdot , \theta^*(y) \rangle_{\mathcal H} \|  \leq
      \| \theta \| \| y \|
    \]
for any element $y \in l_2(A)$ which belongs to the dense in $l_2(A)$ domain
of $\theta^*$. So the norm of $\theta^*$ is bounded by the same constant as
the norm of $\theta$, and $\theta^*$ can be considered as a bounded
$A$-linear map of ${\mathcal H}$ into ${\mathcal H}'$.

Applying $\theta^*$ to the dense in $l_2(A)$ subset of all finite $A$-linear
combinations of the elements $\{ e_j : j \in \mathbb J \}$ the corresponding
range can be seen to be contained in the standard copy of ${\mathcal H}$
inside ${\mathcal H}'$. Hence, the entire image of $\theta^*$
has to belong to the norm-closed set ${\mathcal H} \hookrightarrow {\mathcal
H}'$. This shows the correctness of the definition and the existence of
$\theta^*$ as an adjoint operator of $\theta$. Finally, because $\theta$ is
adjointable, injective and has closed range the operator $\theta^*$ is
surjective, cf.~\cite[Th.~15.3.8]{Wegge-Olsen}.

\smallskip
Since the operator $\theta$ is now shown to be adjointable, injective, bounded
from below and admitting a closed range, the Hilbert $A$-module $l_2(A)$ splits
into the orthogonal sum $l_2(A)= \theta({\mathcal H}) \oplus {\rm Ker}(\theta^*)$
by \cite[Th.~15.3.8]{Wegge-Olsen}. Denote the resulting orthogonal projection
of $l_2(A)$ onto $\theta({\mathcal H})$ by $P$. We want to show that $P(e_j)=
\theta(x_j)$ for every $j \in {\mathbb J}$. For every $x \in \mathcal H$ the
following equality is valid:
 \begin{eqnarray} \label{eqn-adjoint}
   \langle \theta(x),P(e_j) \rangle_{\theta({\mathcal H})} & = &
        \langle P(\theta(x)),e_j \rangle_{l_2} \\ \nonumber
      & = & \langle \theta(x),e_j \rangle_{l_2} \\ \nonumber
      & = & \langle x,x_j \rangle_{{\mathcal H}} \\ \nonumber
      & = & \langle \theta(x),\theta(x_j) \rangle_{\theta({\mathcal H})} \, .\\
      \nonumber
 \end{eqnarray}
In the third equality of the equation above the fact was used that
$\langle \theta(y),e_j \rangle_{l_2} = \langle y,x_j \rangle_{{\mathcal H}}$
for every $y \in \mathcal H$ by the definition of $\theta$. Since $(P(e_j)-
\theta(x_j)) \in \theta({\mathcal H})$ and $x \in \mathcal H$ is arbitrarily
chosen the identity $P(e_j)=\theta(x_j)$ follows for every $j \in {\mathbb J}$.

\smallskip
Since $\theta({\mathcal H})$ is generated by the set $\{ \theta(x_j) : j \in
\mathbb J \}$ and since $\theta$ is an isometry the Hilbert $A$-module
$\mathcal H$ is generated by the set $\{ x_j : j \in \mathbb J \}$ as a Banach
$A$-module. By \cite[Example ${\rm A}_1$]{HanLarson} a standard normalized
tight frame in a finite-dimensional Hilbert space $H$ can contain more
non-zero elements than the dimension of $H$. So the zero element of $\mathcal
H$ may admit a non-trivial decomposition $0= \sum_j a_jx_j$ for some elements
$\{ a_j : j \in \mathbb J \} \subset A$ in some situations.
\end{proof}

\begin{corollary} \label{on-basis}
  Let $A$ be a unital C*-algebra, $\{ {\mathcal H}, \langle .,. \rangle \}$ be a
  finitely or countably generated Hilbert $A$-module.
  Suppose that $\{ x_j : j \in {\mathbb J} \}$ is a standard Riesz basis for
  $\mathcal H$ that is a normalized tight frame. Then $\{ x_j : j \in
  {\mathbb J} \}$ is an orthogonal Hilbert basis with the additional property
  that $\langle x_j,x_j \rangle = \langle x_j,x_j \rangle^2$ for any $j \in
  \mathbb J$.  The converse assertions holds too.
\end{corollary}

\begin{proof}
Since $\{ x_j : j \in {\mathbb J} \}$ is a normalized tight frame we get
$x_j = \sum_i \langle x_j,x_i \rangle x_i$ for any $j \in \mathbb J$ by the
reconstruction formula. The basis property forces $\langle x_j,x_i \rangle x_i
= 0$ for any $i \not= j$ and each fixed $j$. However, the right carrier projection
of $\langle x_j,x_i \rangle$ equals the carrier projection of $x_i$ for every
$i \in \mathbb J$ if calculated inside the bidual von Neumann algebra $A^{**}$.
So $\langle x_j,x_i \rangle =0$ for any $i \not= j$. Proposition \ref{on-basis2}
gives the second property of the Hilbert basis. The converse implication is a
simple calculation fixing an element $x \in \mathcal H$ and setting $x =
\sum_j a_jx_j$ for some elements $\{ a_j : j \in \mathbb J \} \subset A$ and
the given orthonormal basis $\{ x_j : j \in \mathbb J \}$ of $\mathcal H$:
   \begin{eqnarray*} 
       \langle x,x \rangle & = &
               \left\langle \sum_{j \in {\mathbb J}} a_j x_j ,
               \sum_{k \in {\mathbb J}} a_k x_k  \right\rangle 
           =   \sum_{j \in {\mathbb J}} a_j \langle x_j,x_j \rangle a_j^* \\
         & = & \sum_{j \in {\mathbb J}} a_j \langle x_j,x_j \rangle^2 a_j^*
           =   \sum_{j \in {\mathbb J}} \langle a_jx_j,x_j \rangle \langle x_j,
               a_jx_j \rangle \\
         & = & \sum_{j \in {\mathbb J}} \left\langle \sum_{k \in {\mathbb J}}
               a_k x_k , x_j \right\rangle \left\langle x_j ,
               \sum_{l \in {\mathbb J}} a_l x_l  \right\rangle 
           =   \sum_{j \in {\mathbb J}} \langle x,x_j \rangle \langle x_j,x
               \rangle \\       
   \end{eqnarray*}
Note that we applied the supposed equality $\langle x_j, x_j \rangle =
\langle x_j, x_j \rangle^2$, $ j \in \mathbb J$, as the third transformation
step. Since $x \in \mathcal H$ is arbitrarily selected the special orthogonal
basis $\{ x_j : j \in \mathbb J \}$ turns out to be a normalized tight frame
and hence, a Riesz basis.
\end{proof}

We have an easy proof of the uniqueness of the $A$-valued inner product with
respect to which a given frame is normalized tight, generalizing a fact known
for orthonormal Hilbert bases. Note that standard frames can be replaced by
general frames in Corollary \ref{on-basis} as additional investigations
show at section eight.

\begin{corollary}
   Let $A$ be a unital C*-algebra, $\mathcal H$ be a finitely or countably
   generated Hilbert $A$-module, and $\{ x_j : j \in \mathbb J \}$ be a
   standard frame of $\mathcal H$.
   Assume that this frame is normalized tight with respect to two $A$-valued
   inner products $\langle .,. \rangle_1$, $\langle .,. \rangle_2$ on
   $\mathcal H$ that induce equivalent norms to the given one. Then $\langle
   x,y \rangle_1 = \langle x,y \rangle_2$ for any $x,y \in \mathcal H$.
   In other words, the $A$-valued inner product with respect to which a
   standard frame is normalized tight is unique.
\end{corollary}

\begin{proof}
By supposition and Theorem \ref{th-reconstr} we have the reconstruction
formulae
\begin{equation} \label{equ-111}
   x = \sum_{j \in \mathbb J} \langle x,x_j \rangle_1 x_j \, , \,
   y = \sum_{j \in \mathbb J} \langle y,x_j \rangle_2 x_j
\end{equation}
for any $x,y \in \mathcal H$. Taking the $A$-valued inner product of $x$ by
$y$ with respect to $\langle .,. \rangle_2$ and the $A$-valued inner product
of $y$ by $x$ with respect to $\langle .,. \rangle_1$ simultaneously the
right sides of (\ref{equ-111}) become adjoint to one another elements of $A$.
Since $x,y$ are arbitrarily selected elements of $\mathcal H$ the coincidence
of the inner products is demonstrated.
\end{proof}

Remarkably the frame transform of any standard frame preserves the crucial
operator properties known for frame transforms of Hilbert space theory.

\begin{theorem}  \label{th-orthogonal} {\rm (frame transform) } \newline
  Let $A$ be a unital C*-algebra, $\{ {\mathcal H}, \langle .,. \rangle \}$ be a
  finitely or countably generated Hilbert $A$-module.
  Suppose that $\{ x_j : j \in {\mathbb J} \}$ is a standard frame for $\mathcal H$.
  Then the corresponding frame transform $\theta: {\mathcal H} \to l_2(A)$
  defined by $\theta(x)=\{ \langle x,x_j \rangle \}_{j \in {\mathbb J}}$
  $(x \in {\mathcal H})$ possesses an adjoint operator and realizes an
  embedding of $\mathcal H$ onto an orthogonal summand of $l_2(A)$.
  The formula $\theta^*(e_j)=x_j$ holds for every $j \in {\mathbb J}$.
\end{theorem}

\begin{proof}
The set $\{ x_j : j \in {\mathbb J} \}$ is supposed to be standard frame for
the Hilbert $A$-module $\mathcal H$. Refering to the definition of module
frames we have the inequality
\[
   C \cdot \langle x,x \rangle \leq \sum_{j \in \mathbb J} \langle x,x_j \rangle
   \langle x_j,x \rangle = \langle \theta(x),\theta(x) \rangle \leq D \cdot
   \langle x,x \rangle
\]
valid for every $x \in \mathcal H$ and two fixed numbers $0 < C,D$. So the
image of $\theta$ inside $l_2(A)$ has to be closed since $\mathcal H$ is
closed by assumption and the operator $\theta$ is bounded from above and
from below.

The proof of the existence of an adjoint to $\theta$ operator $\theta^* :
l_2(A) \to \mathcal H$ is exactly the same as given in the case of normalized
tight frames, cf.~proof of Theorem \ref{th-reconstr}. Also, the arguments for
$\theta({\mathcal H})$ being an orthogonal summand of $l_2(A)$ can be repeated
as given there.
\end{proof}

For an extended reconstruction formula we refer to Theorem \ref{prop-canondual}
below since some more investigations are necessary to establish it.

\begin{corollary} \label{basis-sum}
{\rm (cf.~\cite[Prop.~2.8]{Casazza:1}) } \newline
If $\{ x_j : j \in {\mathbb J} \}$ is a standard normalized tight frame in a
Hilbert $A$-module $\mathcal H$ then $\{ \theta(x_j) : j \in {\mathbb J} \}$
is the average of two orthonormal Hilbert bases of the Hilbert $A$-module
$l_2(A)$.  \newline
More precisely, let $\{ e_j : j \in \mathbb J \}$ be a fixed Riesz basis of
$\mathcal H$ and at the same time a standard normalized tight frame. Then 
$\theta(x_j) = 1/2 \cdot [(P(e_j)+(1-P)(e_j)) + (P(e_j) - (1-P)(e_j))]$ for
every $j \in {\mathbb J}$ and $P: l_2(A) \to \theta({\mathcal H})$ the
respective orthogonal projection.
\end{corollary}

Since the short proof is straightforward we only mention that $(2P-1)$
is a self-adjoint isometry forcing $\{ (2P-1)(e_j) : j \in \mathbb J \}$ to
be a Riesz basis of the same kind as $\{ e_j : j \in \mathbb J \}$.

\begin{corollary}  \label{cor-43}
  Let $\{ x_j : j \in {\mathbb J} \}$ be an orthogonal Hilbert basis of a
  finitely or countably generated Hilbert $A$-module $\mathcal H$ with
  the property $\langle x_j,x_j \rangle = \langle x_j,x_j \rangle^2$. For
  every partial isometry $V \in {\rm End}_A^*({\mathcal H})$ the sequence
  $\{ V(x_j) : j \in {\mathbb J} \}$ becomes a standard normalized tight frame of
  $V({\mathcal H})$.
\end{corollary}

\begin{proof}
Since $\{ x_j : j \in {\mathbb J} \}$ is an orthogonal Hilbert basis
of $\mathcal H$ with $\langle x_j,x_j \rangle = \langle x_j,x_j \rangle^2$
$\{ x_j : j \in {\mathbb J} \}$ has the property of a standard normalized
tight frame. Writing down this property for the special setting $x=V^*V(y)$
we obtain
\begin{eqnarray*}
\sum_n \langle V(y),V(e_j) \rangle \langle V(e_j),V(y) \rangle & = &
   \sum_n \langle V^*V(y),e_j \rangle \langle e_j,V^*V(y) \rangle \\
   & = & \langle V^*V(y),V^*V(y) \rangle \\
   & = & \langle V^*V(y),y \rangle \\
   & = & \langle V(y),V(y) \rangle
\end{eqnarray*}
\end{proof}

For normalized tight frames $\{ y_j : j \in \mathbb J \}$ in finite-dimensional
Hilbert spaces $H$ we have a ``magic'' formula: $\sum_j \langle y_j,y_j
\rangle = {\rm dim}(H)$, without further requirements to the frame,
cf.~\cite[Cor.~1.2, (iii)]{HanLarson}. Example \ref{ex5} tells us that we
cannot expect a full analogy of this fact for finitely generated
Hilbert C*-modules over non-commutative C*-algebras. But, the formula also
does not survive in a weak sense, for example, giving the same sum value for
\linebreak[4]
every frame with the same number of non-zero elements, cf.~Example \ref{ex5}
and a frame        \linebreak[4]
$\{ 1_A \otimes \sqrt{2}^{-1},1_A \otimes \sqrt{2}^{-1} \}$ for
$A={\rm B}(l_2)$.
However, if the underlying C*-algebra is commutative a similar ''magic''
formula can still be obtained.

\begin{proposition} \label{prop-magic} {\rm (the ''magic'' formula)} \newline
   Let $A={\rm C}(X)$ be a commutative unital C*-algebra, where $X$ is the
   appropriate compact Hausdorff space.
   For any finitely generated Hilbert $A$-module $\mathcal H$ and any standard
   normalized tight frame $\{ y_j : j \in \mathbb J \}$ of $\mathcal H$ the
   (weakly converging) sum $\sum_j \langle y_j,y_j \rangle$ results in a
   continuous function on $X$ with constant non-negative integer values on
   closed-open subsets of $X$. The limit does not depend on the choice of the
   normalized tight frame of $\mathcal H$.
\end{proposition}

\begin{proof}
Consider a normalized tight frame $\{ z_j : j \in \mathbb J \}$ of $\mathcal
H$. For this normalized tight frame the sum exists as a weak limit in $A^{**}$.
Fixing a point $x_0 \in X$ and applying the Hilbert space formula to the
Hilbert space frame $\{ z_j(x_0) : j \in \mathbb J \}$ we obtain 
$\sum_j \langle z_j(x_0),z_j(x_0) \rangle \in \mathbb N$,
\cite[Cor.~1.2, (iii)]{HanLarson}. Therefore, the sum is locally constant
because the number obtained is precisely the dimension of the fibre over $x_0$
in the dual to $\mathcal H$ locally trivial vector bundle over $X$, and the
dimension of fibres is locally constant (cf.\cite[\S 13]{Wegge-Olsen}).
Since closed-open subsets of $X$ are compact we get the desired properties of
the resulting function on $X$ in this particular case.

For an arbitrary standard normalized tight frame $\{ y_j : j \in \mathbb J \}$
for $\mathcal H$ we can again fix a point $x_0 \in X$. Comparing the sums
$\sum_j \langle z_j(x_0),z_j(x_0) \rangle$ and $\sum_j \langle y_j(x_0),
y_j(x_0) \rangle$ we obtain their equality by \cite[Cor.~1.2, (iii)]{HanLarson}. Since
$x_0 \in X$ was arbitrarily chosen the statement follows.
\end{proof}

To understand this ''magic'' formula
This ``magic'' formula is similar to the dimension formula for frames in
finite-dimensional Hilbert spaces (cf.~\cite[Cor.~1.2, (iii)]{HanLarson}).
To understand the formula we had to use the categorical equivalence
between locally trivial vector bundles over $X$ and finitely generated Hilbert
${\rm C}(X)$-modules known as Serre-Swan's theorem \cite{Serre,Swan}.
Interpretating $\mathcal H$ as a set of continuous sections of a locally trivial
vector bundle over $X$ the formula describes the dimension
of the fibre over every point of the base space $X$ in this vector bundle.
Unfortunately, the lack of a localization principle in the non-commutative
case does not allow to find analogous formulae for frames of finitely generated
Hilbert $A$-modules over non-commutative C*-algebras $A$.

\smallskip
Another field of applications of frames are Hilbert-Schmidt operators on
finitely or countably generated Hilbert $A$-modules $\mathcal H$ over unital
commutative C*-algebras $A$ (cf.~\cite{DH}). Since $\mathcal H$ contains a
standard normalized tight frame $\{ x_j : j \in \mathbb J \}$ by Kasparov's
theorem \cite[Th.~1]{Kas} and Corollary \ref{cor-43} we can say the following:
an adjointable bounded $A$-linear operator $T$ on $\mathcal H$ is {\it (weakly)
Hilbert-Schmidt} if the sum $\sum_{j} \langle T(x_j),T(x_j)
\rangle$ converges weakly. This definition is justified by the following
fact:

\begin{proposition}   \label{prop-HS}
  Let $A$ be a unital commutative C*-algebra, $\mathcal H$ be a finitely or
  countably generated Hilbert $A$-module, and $\{ x_j  : j \in \mathbb J \}$
  and $\{ y_j : j \in \mathbb J \}$ be two standard normalized tight frames of
  $\mathcal H$. Consider an adjointable bounded $A$-linear operator $T$ on
  $\mathcal H$.
  If the sum $\sum_{j} \langle T(x_j),T(x_j) \rangle$ converges
  weakly then the sum $\sum_{j} \langle T(y_j), T(y_j) \rangle$
  also converges weakly and gives the same value in $A^{**}$. Furthermore, if
  $T$ is replaced by $T^*$ then the value of this sum does not change.
\end{proposition}

\begin{proof}
We have only to check a chain of equalities in $A^{**}$ that is valid for our
standard normalized tight frames. For an arbitrary fixed standard normalized
tight frame $\{ z_k : k \in \mathbb J \}$ we have
  \begin{eqnarray*}
    \sum_{j} \langle T(x_j),T(x_j) \rangle & = &
      \sum_{k} \sum_{j}
      \langle T(x_j),z_k \rangle \langle z_k,T(x_j) \rangle
     =
      \sum_{k} \sum_{j}
      \langle x_j,T^*(z_k) \rangle \langle T^*(z_k),x_j \rangle \\
    & = &
      \sum_{j} \sum_{k}
      \langle T^*(z_k),x_j \rangle \langle x_j,T^*(z_k) \rangle
     =
      \sum_{k} \langle T^*(z_k),T^*(z_k) \rangle
  \end{eqnarray*}
in case one of the sums at either ends converges weakly. Since we can repeat
our calculations for the other standard normalized tight frame $\{ y_j : j
\in \mathbb J \}$ and since we can choose $z_j=x_j$ for all $j \in \mathbb J$ the
statement of the proposition follows.
\end{proof}

This proposition might be new even for Hilbert spaces and for the definition
of the Hilbert-Schmidt norm of Hilbert-Schmidt operators there. Unfortunately,
the commutativity of the C*-algebra $A$ cannot be omitted.

%%%%%%%%%%%%%%%%%%%%%%%%%%%%%%%%%%%%%%%%%%%%%%%%%%%%%%%%%%%%%%%%%%%%%
\section{Complementary frames, unitary equivalence and similarity}
\label{sec-complement}

In this section we consider geometrical dilation results for frames in
Hilbert C*-modules. The central two concepts are: $\,$ (i) the inner direct
sum of frames with respect to a suitable embedding of the original Hilbert
C*-module into a larger one as an orthogonal summand and $\,$ (ii) the
existence of a complementary frame in the orthogonal complement of this
embedding.
The description of the Hilbert space results can be found in \cite{HanLarson}
as Corollary 1.3, Propositions 1.4-1.7 and 1.9. A more detailed account to
inner sum decompositions of module frames can be found in \cite{FL99}.

\begin{proposition}  \label{prop-complement}
  Let $A$ be a unital C*-algebra, $\mathcal H$ be a finitely (resp., countably)
  generated Hilbert $A$-module and $\{ x_j : j \in {\mathbb J} \}$ be a
  standard normalized tight frame in $\mathcal H$.
  Then there exists another countably generated Hilbert $A$-module $\mathcal M$
  and a standard normalized tight frame $\{ y_j : j \in {\mathbb J} \}$ in
  $\mathcal M$ such that the sequence
  \[
     \{ x_j \oplus y_j : j \in {\mathbb J} \}
  \]
  is an orthogonal Hilbert basis for the countably generated Hilbert
  $A$-module ${\mathcal H} \oplus {\mathcal M}$ with the property
  $\langle x_j \oplus y_j, x_j \oplus y_j \rangle = \langle x_j \oplus y_j,
  x_j \oplus y_j \rangle^2$ for every $j \in \mathbb J$. The complement
  $\mathcal M$ can be selected in such a way that ${\mathcal H} \oplus
  {\mathcal M} = l_2(A)$ and hence, $1_A = \langle x_j \oplus y_j, x_j \oplus y_j
  \rangle$.
  \newline
  If $\mathcal H$ is finitely generated and the index set ${\mathbb J}$ is finite
  then $\mathcal M$ can be chosen to be finitely generated, too, and ${\mathcal H}
  \oplus {\mathcal M} = A^N$ for $N = |{\mathbb J}|$.
  \newline
  If $\{ x_j : j \in \mathbb J \}$ is already an orthonormal basis then
  ${\mathcal M}= \{ 0 \}$, i.e.~no addition to the frame is needed.
  If ${\mathbb J}$ is finite and $\mathcal M$ is not finitely generated then
  infinitely many times $0_{\mathcal H}$ has to be added to the frame $\{ x_j :
  j \in \mathbb J \}$ to make sense of the statement.
\end{proposition}

\begin{proof}
By Theorem \ref{th-reconstr} there is a standard isometric embedding
of $\mathcal H$ into $l_2(A)$ induced by the frame transform $\theta$.
In the context of that embedding
$\theta({\mathcal H})$ is an orthogonal summand of $l_2(A)$, and the
$A$-valued inner products on $\mathcal H$ and on $\theta({\mathcal H})$
coincide.
The corresponding projection $P: l_2(A) \to \theta({\mathcal H})$ maps
the standard orthonormal Riesz basis $\{ e_j : j \in {\mathbb J} \}$ of
$l_2(A)$ onto the frame $\{ \theta(x_j) : j \in {\mathbb J} \}$. Set
${\mathcal M} = (I-P)(l_2(A))$ and consider $y_j=(I-P)(e_j)$ for $j \in
{\mathbb J}$. These objects possess the required properties.

If $|{\mathbb J}|$ is finite the frame transform $\theta$ can take its image
in the standard Hilbert $A$-submodule $A^N \subset l_2(A)$ with $N=|{\mathbb J}|$.
\end{proof}

\begin{proposition} \label{prop-unitary-equiv}
  Let $A$ be a unital C*-algebra, $\mathcal H$ be a countably generated Hilbert
  $A$-module and $\{ x_j : j \in {\mathbb J} \}$ be a standard normalized tight
  frame for $\mathcal H$, where the index set ${\mathbb J}$ is countable or finite.
  Suppose, there exist two countably generated Hilbert $A$-modules $\mathcal M$,
  $\mathcal N$ and two normalized tight frames $\{ y_j : j \in {\mathbb J} \}$,
  $\{ z_j : j \in {\mathbb J} \}$ for them, respectively, such that
  \[
     \{ x_j \oplus y_j : j \in {\mathbb J} \} \, ,
     \, \{ x_j \oplus z_j : j \in {\mathbb J} \}
  \]
  are orthogonal Hilbert bases for the countably generated
  Hilbert $A$-modules ${\mathcal H} \oplus {\mathcal M}$, ${\mathcal H} \oplus
  {\mathcal N}$, respectively, where we have the value properties
  $\langle x_j \oplus y_j, x_j \oplus y_j \rangle = \langle x_j \oplus y_j,
  x_j \oplus y_j \rangle^2$ and $\langle x_j \oplus z_j, x_j \oplus z_j \rangle
   = \langle x_j \oplus z_j, x_j \oplus z_j \rangle^2$. If $\langle y_j,y_j
   \rangle_{{\mathcal M}} = \langle z_j,z_j \rangle_{{\mathcal N}}$
  for every $j \in {\mathbb J}$, then there exists a unitary transformation
  $U: {\mathcal H} \oplus {\mathcal M} \to {\mathcal H} \oplus {\mathcal N}$
  mapping $\mathcal M$ onto $\mathcal N$ and satisfying $U(y_j)=z_j$ for every
  $j \in {\mathbb J}$.

  The additional remarks of Proposition \ref{prop-complement} apply in the
  situation of finitely generated Hilbert $A$-modules appropriately.
\end{proposition}

\begin{proof}
Set $e_j = x_j \oplus y_j$ and $f_j = x_j \oplus z_j$ and
define $U'(e_j)=f_j$. By assumption the $A$-valued inner products are
preserved by $U'$, and $U'$ extends to a unitary map between ${\mathcal H}
\oplus {\mathcal M}$ and ${\mathcal H} \oplus {\mathcal N}$ by $A$-linearity.
Fix $x \in \mathcal H$. Then the equality
\[
  \langle x,x_j \rangle_{{\mathcal H}}=\langle x \oplus 0_{{\mathcal M}}, e_j
  \rangle = \langle x \oplus 0_{{\mathcal N}}, f_j \rangle  \, , \, \, j \in
  {\mathbb J} \, ,
\]
is valid. So $x \oplus 0_{{\mathcal M}}= \sum_j \langle x \oplus
0_{{\mathcal M}}, e_j \rangle e_j = \sum_j \langle x,x_j \rangle e_j$ and
$x \oplus 0_{{\mathcal N}} = \sum_j \langle x,x_j \rangle f_j$ for $j \in
{\mathbb J}$. Applying $U'$ the equality $U'(x \oplus 0_{{\mathcal M}})=x
\oplus 0_{{\mathcal N}}$ yields.
Consequently, $U'$ splits into the direct sum of the identity mapping on
the first component and of a unitary operator $U: {\mathcal M} \to
{\mathcal N}$ on the second component.
\end{proof}

\begin{theorem}  \label{prop-dualframe1}
  Let $\{ x_j : j \in {\mathbb J} \}$ be a standard frame of a finitely or
  countably generated Hilbert $A$-module $\mathcal H$. Then $\{ x_j : j \in
  {\mathbb J} \}$ is the image of a standard normalized tight frame $\{ y_j :
  j \in {\mathbb J} \}$ of $\mathcal H$ under an invertible adjointable bounded
  $A$-linear operator $T$ on ${\mathcal H}$. The operator $T$ can be chosen to
  be positive and equal to the square root of $\theta^*\theta$, where
  $\theta$ is the frame transform corresponding to $\{ x_j \}$.
  \newline
  Conversely, the image of a standard normalized tight frame $\{ y_j : j \in
  {\mathbb J} \}$ of $\mathcal H$ under an invertible adjointable bounded $A$-linear
  operator $T$ on ${\mathcal H}$ is a standard frame of $\mathcal H$.
  \newline
  The frame $\{ x_j \}$ is a set of generators of $\mathcal H$ as an Hilbert
  $A$-module. The frame elements do not form a Hilbert basis, in general.
\end{theorem}

\begin{proof}
If $T$ is an invertible adjointable bounded $A$-linear operator
on $\mathcal H$ and $\{ y_j : j \in {\mathbb J} \}$ is a standard normalized tight
frame of $\mathcal H$, then the sequence $\{ x_j=T(y_j) : j \in {\mathbb J} \}$
fulfills the equality
 \begin{eqnarray}       \label{eq3}
   \sum_j \langle x,x_j \rangle \langle x_j,x \rangle & = &
   \sum_j \langle x,T(y_j) \rangle \langle T(y_j),x \rangle  \\ \nonumber
   & = &  \sum_j \langle T^*(x),y_j \rangle \langle y_j,T^*(x) \rangle \\ \nonumber
   & = & \langle T^*(x),T^*(x) \rangle \\ \nonumber
 \end{eqnarray}
for every $x \in \mathcal H$. Since $\|T^{-1}\|^{-2} \langle x,x \rangle \leq
\langle T^*(x),T^*(x) \rangle \leq \|T\|^2 \langle x,x \rangle$ for every
$x \in \mathcal H$ (cf.~\cite{Pa1}) and since the sum in (\ref{eq3}) converges
in norm, the sequence $\{ x_j : j \in {\mathbb J} \}$ is a standard frame of
$\mathcal H$ with frame bounds $C \geq \|T^{-1}\|^{-2}$ and $D \leq \|T\|^2$.

Conversely, for an arbitrary standard frame $\{ x_j : j \in {\mathbb J} \}$ of a
countably generated Hilbert $A$-module $\mathcal H$ the frame transform
$\theta: {\mathcal H} \to l_2(A)$, $\theta(x) = \{ \langle x,x_j \rangle : j
\in {\mathbb J} \}$, is adjointable by Theorem \ref{th-orthogonal}. Moreover,
$\theta^*$ restricted to the orthogonal summand $\theta({\mathcal H})$ of $l_2(A)$
is an invertible operator as $\theta^*$ is the adjoint operator of $\theta$,
where $\theta$ has to be regarded
as an invertible operator from $\mathcal H$ to $\theta({\mathcal H})$.
So the mapping $\theta^*\theta$ becomes an invertible
positive bounded $A$-linear operator onto $\mathcal H$, and the equality
\[
\langle \theta(x),\theta(x) \rangle_{l_2} =
   \sum _j \langle x,x_j \rangle_{{\mathcal H}} \langle x_j,x \rangle_{{\mathcal H}}
\]
holds for every $x \in \mathcal H$. Set $y_x = (\theta^*\theta)^{1/2}(x)$ for
each $x \in \mathcal H$, $y_j = (\theta^*\theta)^{-1/2}(x_j)$ for $j \in {\mathbb J}$.
Then the equality
\[
   \langle y_x,y_x \rangle_{{\mathcal H}} = \langle \theta(x),\theta(x)
      \rangle_{l_2}
   =  \sum _j \langle x,x_j \rangle_{{\mathcal H}} \langle x_j,x
   \rangle_{{\mathcal H}}
   =  \sum_j \langle y_x, y_j \rangle_{{\mathcal H}} \langle y_j,y_x
   \rangle_{{\mathcal H}}
\]
is valid since $x \in \mathcal H$ was arbitrarily chosen and the sum on the
right side converges in norm by supposition.
So the sequence $\{ y_j=(\theta^*\theta)^{-1/2}(x_j) : j \in {\mathbb J} \}$
has been characterized as a standard normalized tight frame of $\mathcal H$.
The operator $T=(\theta^*\theta)^{1/2}$ is the
sought operator mapping the standard normalized frame $\{ y_j \}$ onto
the standard frame $\{ x_j \}$.

The property of a standard frame to be a set of generators for $\mathcal H$ as
a Hilbert $A$-module can be derived from the analogous property of standard
normalized tight frames which is preserved under adjointable invertible
mappings, cf.~Theorem \ref{th-reconstr}.
\end{proof}

\begin{remark} {\rm
  Applying the techniques described in the appendix, we can show that
  the image of a standard normalized tight frame under a non-adjointable invertible
  bounded $A$-linear operator $T$ on $\mathcal H$ is still a frame of $\mathcal H$
  with $C \geq \|T^{-1}\|^{-2}$, $D\leq \|T\|^2$. However, the adjoint
  operator $T^*$ needed for calculations exists as an element of the
  W*-algebra ${\rm End}_A^*(({\mathcal H}^\#)')$ only. In other words,
  there exists an element $x \in \mathcal H$ such that the left side sum
  in (\ref{eq3}) does not converge in norm since $T^*(x) \not\in \mathcal H$.
  The resulting frame $\{ x_j = T(y_j) \}$ turns out to be non-standard.}
\end{remark}

\begin{corollary} {\rm (cf.~\cite[Prop.~2.9]{Casazza:1}) } \newline
  Every standard frame in a Hilbert $A$-module $\mathcal H$ is similar to another
  standard frame in $\mathcal H$ which is mapped to the average of two orthonormal
  bases of $l_2(A)$ by its frame transform.
\end{corollary}

For proof arguments we refer to the Theorems \ref{th-reconstr},
\ref{prop-dualframe1} and Corollary \ref{basis-sum}.

\begin{proposition}   \label{prop-dualframe2}
  Let $\{ x_j : j \in {\mathbb J} \}$ be a standard frame of a finitely or
  countably generated Hilbert $A$-module $\mathcal H$. There exists a Hilbert
  $A$-module $\mathcal M$ and a normalized tight frame $\{y_j : j \in {\mathbb J} \}$
  in $\mathcal M$ such that the sequence $\{ x_j \oplus y_j : j \in {\mathbb J} \}$
  is a standard Riesz basis in ${\mathcal H} \oplus {\mathcal M}$ with the same frame
  bounds for $\{ x_j \}$ and $\{ x_j \oplus y_j \}$. The Hilbert $A$-module
  $\mathcal M$ can be chosen in such a way that ${\mathcal H} \oplus {\mathcal M} =
  l_2(A)$. If $\mathcal H$ is finitely generated and the index set ${\mathbb J}$ is
  finite then $\mathcal M$ can be chosen to be finitely generated, too, and
  ${\mathcal H} \oplus {\mathcal M} = A^N$ for $N=|{\mathbb J}|$.
  \newline
  In general, $\mathcal M$ cannot be chosen as a submodule of $\mathcal H$, and the
  resulting standard Riesz basis may be non-orthogonal. A uniqueness result
  like that one in Proposition \ref{prop-unitary-equiv} fails to be true, in
  general.
\end{proposition}

\begin{proof}
By Theorem \ref{prop-dualframe1} there exists a standard normalized tight
frame $\{ z_j : j \in \mathbb J \}$ for $\mathcal H$ and an adjointable
invertible operator $T$ on $\mathcal H$ such that $x_j=T(z_j)$ for any $j \in
{\mathbb J}$.
Moreover, there is another Hilbert $A$-module $\mathcal M$ and a standard
normalized tight frame $\{ y_j : j \in \mathbb J \}$ for $\mathcal M$ such
that the sequence $\{z_j \oplus y_j : j \in \mathbb J \}$ is an orthogonal
Hilbert basis in ${\mathcal H} \oplus {\mathcal M}$, see Proposition
\ref{prop-complement}. Then $T \oplus {\rm id}$ is an adjointable invertible
operator on ${\mathcal H} \oplus {\mathcal M}$ mapping the sequence $\{ z_j
\oplus y_j : j \in \mathbb J \}$ onto the sequence $\{x_j \oplus y_j :
j \in \mathbb J\}$.
Hence, the latter is a standard Riesz basis for ${\mathcal H} \oplus {\mathcal M}$
according to Theorem \ref{prop-dualframe1}. The statement regarding bounds
is obvious, the special choices for $\mathcal M$ can be derived from the
reconstruction formula. The additional remarks have been already shown to be
true for particular Hilbert space situations in
\cite[Prop.~1.6, Example B]{HanLarson}.
\end{proof}

\begin{corollary}   \label{cor-Riesz}
  Let $\{ x_j : j \in {\mathbb J} \}$ be a standard Riesz basis of a finitely or
  countably generated Hilbert $A$-module $\mathcal H$. Then $\{ x_j : j \in
  {\mathbb J} \}$ is the image of a standard normalized tight frame and
  Hilbert basis $\{ y_j : j \in {\mathbb J} \}$
  of $\mathcal H$ under an invertible adjointable bounded $A$-linear operator $T$
  on ${\mathcal H}$, i.e.~of an orthogonal Hilbert basis $\{ y_j : j \in {\mathbb
  J} \}$ with the property $\langle y_j,y_j \rangle = \langle y_j,y_j \rangle^2$
  for any $j \in \mathbb J$.
  \newline
  Conversely, the image of a standard normalized tight frame and Hilbert basis
  $\{ y_j : j \in {\mathbb J} \}$ of $\mathcal H$ under an invertible adjointable
  bounded $A$-linear operator $T$ on ${\mathcal H}$ is a standard Riesz basis of
  $\mathcal H$.
  \newline
  If a Hilbert $A$-module $\mathcal H$ contains a standard Riesz basis then 
  $\mathcal H$ contains an orthogonal Hilbert basis $\{ x_j : j \in \mathbb J
  \}$ with the frame property $x = \sum_j \langle x,x_j \rangle x_j$ for every
  element $x \in \mathcal H$.
\end{corollary}

Let ${\mathcal H}_1$ and ${\mathcal H}_2$ be Hilbert C*-modules over a fixed
C*-algebra $A$. Let $\{ x_j : j \in \mathbb J \}$ and $\{ y_j : j \in
\mathbb J \}$ be frames for these Hilbert C*-modules, respectively, where the
possibility to select the same index set $\mathbb J$ is essential for our
purposes in the sequel. We call the sequence $\{ x_j \oplus y_j : j \in \mathbb
J \}$ of the Hilbert $A$-module ${\mathcal H}_1 \oplus {\mathcal H}_2$ the
{\it inner direct sum of the frames $\{ x_j : j \in \mathbb J \}$ and $\{ y_j
: j \in \mathbb J \}$}. The two components-frames $\{ x_j : j \in \mathbb J
\}$ and $\{ y_j : j \in \mathbb J \}$ are called {\it inner direct summands}
of the sequence $\{ x_j \oplus y_j : j \in \mathbb J \}$, in particular if the
latter is a frame for ${\mathcal H}_1 \oplus {\mathcal H}_2$.
With these denotations we can reformulate a main result of our investigations
in the following way, cf.~\cite[Th.~1.7]{HanLarson}:

\begin{theorem}
  Standard frames are precisely the inner direct summands of standard Riesz
  bases of $A^N$ or $l_2(A)$. Standard normalized tight frames are precisely
  the inner direct summands of orthonormal Hilbert bases of $A^N$ or $l_2(A)$.
\end{theorem}

The problem whether non-standard frames can be realized as inner direct
summands of generalized Riesz bases of certain canonical Hilbert C*-modules,
or not, is still open. The problem is tightly connected to the existence
problem of a well-behaved frame transform for non-standard frames and
corresponding codomain Banach C*-modules.

\medskip
Proposition \ref{prop-dualframe2} has immediate consequences for the
characterization of algebraically generating sets of (algebraically) finitely
generated Hilbert C*-modules over unital C*-algebras as frames. 
Below we give a transparent proof of the fact that finitely generated Hilbert
$A$-modules over unital C*-algebras $A$ are projective $A$-modules.
Usually, this fact can only be derived from Kasparov's stabilization theorem
for countably generated Hilbert $A$-modules,
cf.~\cite[Cor.~15.4.8]{Wegge-Olsen}.
The smallest appearing number $n \in \mathbb N$ for which a given finitely
generated Hilbert $A$-module is embeddable into the Hilbert $A$-module $A^n$
as an orthogonal summand equals the number of elements of the shortest frame
of the considered Hilbert $A$-module. Also, the general validity of the lower
bound inequality in the chain of inequalities below is a fact possibly not
sufficiently recognized before.

\begin{theorem}  \label{prop-afgHm}
  Every algebraically finitely generated Hilbert $A$-module ${\mathcal H}$
  over a unital C*-algebra $A$ is projective, i.e.~an orthogonal summand of
  some free $A$-module $A^n$ for a finite integer $n \in \mathbb N$.
  Furthermore, any algebraically generating set $\{ x_i : i = 1,...,n \}$ of
  $\mathcal H$ is a frame, and the inequality
  \[
    C \cdot \langle x,x \rangle \leq \sum_{i=1}^n \langle x,x_i \rangle
    \langle x_i,x \rangle \leq D \cdot \langle x,x \rangle
  \]
  holds for every element $x \in \mathcal H$ and some constants $0 < C,D <
  +\infty$.  In other words, the positive bounded module operator $\sum_j
  \theta_{x_j,x_j}$ is invertible.
\end{theorem}

\begin{proof}
Consider the operator $F: A^n \to \mathcal H$ defined by $F(e_i)=x_i$ for
$i=1,...,n$ and for an orthonormal basis $\{ e_i \}_{i=1}^n$ of $A^n$. The
operator $F$ is a bounded $A$-linear, surjective and adjointable operator
since $\mathcal H$ is supposed to be algebraically generated by $\{ x_i : i=
1,...,n \}$ and the Hilbert $A$-module $A^n$ is self-dual,
cf.~\cite[Prop.~3.4]{Pa1}. By \cite[Th.~15.3.8]{Wegge-Olsen} the operator
$F^*$ has to be bounded $A$-linear, injective with closed range. Furthermore,
$F$ possesses a polar decomposition $F=V |F|$, where $A^n= {\rm ker}(F) \oplus
F^*({\mathcal H})$, ${\rm ker}(V)={\rm ker}(F)$ and $V^*({\mathcal H}) =
F^*({\mathcal H})$, see \cite[Th.~15.3.8]{Wegge-Olsen}. The set $\{ V(e_i) :
i=1,...,n \}$ is a normalized tight frame of $\mathcal H$ by Corollary
\ref{cor-43}, and $x_i = (FV^*)(V(e_i))$ for every $i=1,...,n$ by construction.
However, the operator $FV^*$ is invertible on $\mathcal H$. So the set
$\{ x_i : i=1,...,n \}$ is a frame by Proposition \ref{prop-dualframe2}. The
inequality can be obtained from the definition of a frame.
\end{proof}

D.~P.~Blecher pointed out to us that the operator $T=\sum_i \theta_{x_i,x_i}$
is strictly positive by \cite[Cor.~1.1.25]{Jensen/Thomsen}. Since the set
of all ''compact'' module operators on finitely generated Hilbert C*-modules
is a unital C*-algebra $T$ has to be invertible, cf.~\cite[15.O]{Wegge-Olsen}.
This establishes the upper and the lower frame bounds as $\|T\|^2$ and
$\|T^{-1}\|^{-2}$.

We close this subsection with some observations on inner direct sums of
frames. Our interest centers on frame property preserving exchanges of the
second inner direct summand to unitarily equivalent ones.

\begin{proposition}
  Let $A$ be a unital C*-algebra.    \newline
  (i) If $\{ x_j : j \in \mathbb J \}$ is a standard (normalized tight) frame
  for a Hilbert $A$-module $\mathcal H$ and $T$ is a co-isometry on $\mathcal
  H$ (i.e.~$T$ is an adjointable operator such that $T^*$ is an isometry) then
  $\{T(x_j) : j \in \mathbb J \}$ is a standard (normalized tight) frame.
  \newline
  (ii) Let $\{ x_j : j \in \mathbb J \}$ and $\{ y_j : j \in \mathbb J \}$ be
  standard normalized tight frames for Hilbert $A$-modules $\mathcal H$ and
  $\mathcal K$, respectively, that are connected by an adjointable bounded
  operator $T$ obeying the formula $T(x_j)=y_j$ for $j \in {\mathbb J}$. Then
  $T$ is a co-isometry. If $T$ is invertible then it is a unitary.
  \newline
  (iii) Let $\{ x_j : j \in \mathbb J \}$ and $\{ y_j : j \in \mathbb J \}$ be
  standard normalized tight frames for Hilbert $A$-modules $\mathcal H$ and
  ${\mathcal K}$, respectively, with the property that $\{ x_j \oplus y_j :
  j \in \mathbb J\}$ is a standard normalized tight frame for ${\mathcal H}
  \oplus {\mathcal K}$. Then for every standard normalized tight frame
  $\{ z_j : j \in \mathbb J \}$ of the Hilbert $A$-module ${\mathcal K}$ that
  is unitarily equivalent to $\{ y_j : j \in \mathbb J \}$ the sequence
  $\{ x_j \oplus z_j : j \in \mathbb J \}$ again forms a standard normalized
  tight frame of ${\mathcal H} \oplus {\mathcal K}$.
  \newline
  (iv) Let $\{ x_j : j \in \mathbb J \}$ and $\{ y_j : j \in \mathbb J \}$ be
  standard normalized tight frames in Hilbert $A$-modules $\mathcal H$ and
  ${\mathcal K}$, respectively, with the property that $\{ x_j \oplus y_j : j
  \in \mathbb J \}$ is a standard normalized tight frame in ${\mathcal H} \oplus
  {\mathcal K}$. For every standard frame $\{ z_j : j \in \mathbb J \}$
  of the Hilbert $A$-module ${\mathcal K}$ that is similar to $\{ y_j : j \in
  \mathbb J \}$ the sequence $\{ x_j \oplus z_j : j \in \mathbb J \}$ again
  forms a standard frame of ${\mathcal H} \oplus {\mathcal K}$.
\end{proposition}

\begin{proof}
Let $C$ and $D$ be the frame bounds for the standard frame $\{ x_j : j \in
\mathbb J \}$. Then for $x \in \mathcal H$ we obtain the inequality
  \begin{eqnarray*}
    C \cdot \langle x,x \rangle & = & C \cdot \langle T^*(x),T^*(x) \rangle \\
         & \leq &
         \sum_j \langle T^*(x),x_j \rangle \langle x_j,T^*(x) \rangle \\
         & = &
         \sum_j \langle x,T(x_j) \rangle \langle T(x_j),x \rangle \\
         & \leq & D \cdot \langle T^*(x),T^*(x) \rangle = D \cdot \langle x,x \rangle
  \end{eqnarray*}
by E.~C.~Lance's theorem \cite{Lance} and the frame property.
The additional equality in the middle of this chain of two inequalities
introduces a certain expression the comparison of which to both the ends of
the chain of inequalities establishes assertion (i).

Let $\{ x_j : j \in \mathbb J \}$ and $\{ y_j : j \in \mathbb J \}$ be
standard normalized tight frames for Hilbert $A$-modules $\mathcal H$ and
$\mathcal K$, respectively. Suppose there exists an adjointable bounded
operator $T$ such that $T(x_j)=y_j$ for every $j \in {\mathbb J}$. For $y \in
\mathcal K$ the equality
  \[
    \langle T^*(y) T^*(y) \rangle =
    \sum_j  \langle T^*(y),x_j \rangle \langle x_j,T^*(y) \rangle  =
    \sum_j  \langle y,T(x_j) \rangle \langle T(x_j),y \rangle =
    \langle y,y \rangle
  \]
is valid. So $T^*$ is an isometry of the Hilbert $A$-module $\mathcal K$ into
the Hilbert $A$-module $\mathcal H$. If $T$ is invertible then $\mathcal H$
and $\mathcal K$ are unitarily isomorphic by \cite{Lance}. This shows (ii).

To give some argument for (iii) fix a unitary operator $U \in {\rm End}_A({\mathcal K})$ with
the property $U(y_j)=z_j$, $j \in {\mathbb J}$. Then $V={\rm id}
\oplus U \in {\rm End}_A({\mathcal H} \oplus {\mathcal K})$ is a unitary with
the property $V(x_j \oplus y_j)=x_j \oplus z_j$. Hence, the sequence
$\{ x_j \oplus z_j : j \in \mathbb J \}$ is a standard normalized tight
frame for ${\mathcal H} \oplus {\mathcal K}$.
Replacing $U$ by a merely invertible adjointable bounded operator $T$ and
repeating the considerations we obtain assertion (iv).
\end{proof}

%%%%%%%%%%%%%%%%%%%%%%%%%%%%%%%%%%%%%%%%%%%%%%%%%%%%%%%%%%%%%%%%%%%%%%%%%%%%%%%

\section{The canonical dual frame and alternate dual frames}  \label{sec-dual}

The purpose of this section is to establish the existence of canonical and
alternate dual frames of standard frames and to prove fundamental properties
of them. Theorem \ref{prop-canondual} states the general reconstruction
formula for standard frames, the existence of both the frame operator and of
the canonical dual frame. The Propositions 6.2, 6.3, 6.5, 6.6, 6.7 show
relations between canonical dual and alternative dual frames of a given
standard frame.
Example \ref{exam00} below demonstrates one of the differences of
generalized module frame theory for Hilbert C*-modules in comparison to
classical Hilbert space frame theory: the appearance of zero-divisors in most
C*-algebras may cause the non-uniqueness of the dual frame of a standard Riesz
basis.

Let us consider the sequence $\{ (\theta^*\theta)^{-1}(x_j) : j \in {\mathbb
J} \}$ for a standard frame $\{ x_j : j \in {\mathbb J} \}$ for a finitely or
countably generated Hilbert C*-module $\mathcal H$. Denote the map that assigns
to every $x \in \mathcal H$ the corresponding unique pre-image in
$\theta({\mathcal H})$ under $\theta^*$ by $(\theta^*)^{-1}$. This map is
well-defined since $\theta^*$ is injective with image $\mathcal H$. So
$(\theta^*)^{-1}$ is an invertible bounded $A$-linear operator mapping
$\mathcal H$ onto $\theta(\mathcal H)$. Refering to the proof of Theorem
\ref{th-reconstr} and to Theorem \ref{th-orthogonal} we have the following
chain of equalities
  \begin{eqnarray*}
    \theta(x) & = &
      \sum_j \langle \theta(x),e_j \rangle_{l_2} e_j  \\
      & = & \sum_j \langle x, \theta^*(e_j) \rangle e_j
        =   \sum_j \langle x, \theta^*(e_j) \rangle P(e_j) \\
      & = & \sum_j \langle x,x_j \rangle (\theta^*)^{-1}(x_j)
        =   \sum_j \langle \theta(x), (\theta^*)^{-1}(x_j) \rangle_{l_2}
              (\theta^*)^{-1}(x_j)  \\
      & = & \theta \left( \sum_j \langle x,x_j \rangle (\theta^*\theta)^{-1}
              (x_j) \right) 
  \end{eqnarray*}
which holds for every $x \in \mathcal H$ and for the standard orthonormal
Hilbert basis $\{ e_j : j \in \mathbb J \}$ of $l_2(A)$. The pre-last line of
the established equality shows that the sequence $\{(\theta^*)^{-1}(x_j) :
j \in \mathbb J \}$ is a standard normalized tight frame of $\theta({\mathcal
H})$. Since $\theta$ is injective the last line gives a remarkable property
of the sequence $\{ (\theta^*\theta)^{-1}(x_j) : j \in {\mathbb J} \}$:
\[
    x = \sum_j \langle x,x_j \rangle (\theta^*\theta)^{-1} (x_j)
\]
for every $x \in \mathcal H$. Applying $\theta^*$ to this equality and
replacing $x$ by $(\theta^*\theta)^{-1} (x)$ we obtain another equality
dual to the former one:
\[
    x = \sum_j \langle x, (\theta^*\theta)^{-1} (x_j) \rangle x_j
\]
being valid for every $x \in \mathcal H$. We take these two equalities as a
justification to introduce a new notion. 
The frame $\{ (\theta^*\theta)^{-1}(x_j) : j \in {\mathbb J} \}$ is said
to be the {\it canonical dual frame of the frame $\{ x_j : j \in \mathbb J \}$},
and the operator $S=(\theta^*\theta)^{-1}$ is said to be the {\it frame
operator of the frame $\{ x_j : j \in \mathbb J \}$}. In case
the standard frame $\{ x_j : j \in \mathbb J \}$ of $\mathcal H$ is already
normalized tight the operator $S$ is just the identity operator, and the dual
frame coincides with the frame itself.

More generally, we have an existence and uniqueness result (see Theorem below)
that provides us with a reconstruction formula for standard frames. The proof
is only slightly more complicated than in the Hilbert space case
(cf.~\cite[Prop.~1.10, Rem.~1.12]{HanLarson}) since most difficulties were
already overcome establishing the properties of the frame transform.

\begin{theorem}  \label{prop-canondual} {\rm (reconstruction formula)}
  \newline
  Let $\{ x_j : j \in {\mathbb J} \}$ be a standard frame in a finitely or
  countably generated Hilbert $A$-module $\mathcal H$ over a unital C*-algebra $A$.
  Then there exists a unique operator $S \in {\rm End}_A^*({\mathcal H})$ such
  that
  \[
  x = \sum_{j} \langle x,S(x_j) \rangle x_j
  \]
  for every $x \in \mathcal H$. The operator can be explicitely given by the
  formula $S=G^*G$ for any adjointable invertible bounded operator $G$ mapping
  $\mathcal H$ onto some other Hilbert $A$-module $\mathcal K$ and realizing
  $\{ G(x_j) : j \in {\mathbb J} \}$ to be a standard normalized tight frame
  in $\mathcal K$.
  In particular, $S= \theta^{-1} (\theta^*)^{-1} = (\theta^*\theta)^{-1}$ for
  the frame transform $\theta$ with codomain $\theta(\mathcal H)$. So $S$ is
  positive and invertible.
  \newline
  Finally, the canonical dual frame is a standard frame for $\mathcal H$,
  again.
\end{theorem}

\begin{proof}
Let $G \in {\rm End}_A^*({\mathcal H},{\mathcal K})$ be any invertible
operator onto some Hilbert $A$-module $\mathcal K$ with the property that the
sequence $\{ G(x_j) : j \in {\mathbb J} \}$ is a standard normalized tight
frame of $\mathcal K$. The existence of such an operator is guaranteed by
Theorem \ref{th-reconstr} setting ${\mathcal K}=\theta({\mathcal H})$ and
$G=(\theta^*)^{-1}$ (cf.~the introductory considerations of the present
section), or by Theorem \ref{prop-dualframe1}. Set $S=G^*G$ and check the
frame properties of the sequence $\{ S(x_j) : j \in {\mathbb J} \}$:
  \begin{eqnarray*}
    \sum_j \langle x,G^*G(x_j) \rangle x_j & = &
      \sum_j \langle G(x),G(x_j) \rangle x_j =
       \sum_j \langle G(x),G(x_n) \rangle G^{-1}(G(x_j)) \\
    & = & G^{-1} \left( \sum_j \langle G(x),G(x_j) \rangle G(x_j) \right)
      = G^{-1}G(x) = x \, .
  \end{eqnarray*}
The equality implies $\langle S(x),x \rangle = \sum_j \langle x,S(x) \rangle
\langle S(x),x \rangle$ for any $x \in \mathcal H$. Since $G$ is invertible
and $S$ is positive there exist two constants $0 < C,D$ such that the
inequality
\[
   C \cdot \langle x,x \rangle \leq \langle S(x),x \rangle = \sum_j \langle
   x,S(x) \rangle \langle S(x),x \rangle \leq D \cdot \langle x,x \rangle
\]
is fulfilled for every $x \in \mathcal H$. So the sequence $\{ S(x_j) : j \in
{\mathbb J} \}$ is a standard frame of $\mathcal H$ and a dual frame of the
frame $\{ x_j : j \in \mathbb J \}$.

To show the uniqueness of $S$ in ${\rm End}_A^*({\mathcal H})$ and the
coincidence of the found dual frame with the canonical dual frame suppose the
existence of a second operator $T \in {\rm End}_A^*({\mathcal H})$ realizing
the equality $x=\sum_j \langle x,T(x_j) \rangle x_j$ for every $x \in
\mathcal H$. Then we obtain
  \begin{eqnarray*}
    x & = & \sum_j \langle x,T(x_j) \rangle x_j =
              \sum_j \langle x, TG^{-1}G(x_j) \rangle G^{-1}G(x_j) \\
      & = & G^{-1} \left( \sum_j \langle (G^*)^{-1}T^*(x),G(x_j) \rangle
               G(x_j) \right) \\
      & = & G^{-1}((G^*)^{-1}T^*(x)) = (G^*G)^{-1}T^*(x)
  \end{eqnarray*}
for every $x \in \mathcal H$. Consequently, $T=G^*G$ as required.
\end{proof}

If $\{ x_j : j \in \mathbb J \}$ is a standard frame of a Hilbert $A$-module
$\mathcal H$ which is not a Hilbert basis then there may in general exist many
standard frames $\{ y_j : j \in \mathbb J \}$ of $\mathcal H$ for which the
formula
  \begin{equation}  \label{eq-dualframe}
    x= \sum_j \langle x,y_j \rangle x_j
  \end{equation}
is valid. For examples in one- and two-dimensional complex Hilbert
spaces we refer the reader to \cite[\S 1.3]{HanLarson}. We add another example
from C*-theory which reminds the Cuntz algebras ${\mathcal O}_n$:
let $A$ be a C*-algebra with $n$ elements $\{ x_1, ... , x_n \}$ such that
$\sum_i x_i^*x_i =1_A$. Then this set is a standard normalized tight frame of
$A$ by the way of its setting (where $A$ is considered as a left Hilbert
$A$-module). However, any other set $\{ y_1, ... , y_n \}$ of $A$ satisfying
$\sum_i y_i^*x_i = 1_A$ fulfills the analogue of equality (\ref{eq-dualframe})
as well. The choice $y_i=x_i$ is only the one that corresponds to the
canonical dual frame of the initial frame. Other frames can be obtained, for
example, setting $x_1=x_2=\sqrt{2}^{-1} \cdot 1_A$ and $y_1=\sqrt{2} \cdot 1_A$,
$y_2=0_A$. We call the other frames satisfying the equality
(\ref{eq-dualframe}) {\it the alternate dual frames} of a given standard
frame. Note that the frame property of these alternate sequences has to be
supposed since there are examples of non-frame sequences $\{ y_j : j \in
\mathbb J \}$ fulfilling the equality (\ref{eq-dualframe}) in some situations,
\cite[\S 1.3]{HanLarson}.
The following proposition characterizes the operation of taking the canonical
dual frame as an involutive mapping on the set of standard frames,
cf.~\cite[Cor.~1.11]{HanLarson}.

\begin{proposition}
  Let $\{ x_j : j \in \mathbb J \}$ be a standard frame of a Hilbert
  $A$-module $\mathcal H$. Then the canonical dual frame $\{ (\theta^*\theta)^{-1}
  (x_j) : j \in \mathbb J \}$ fulfills the equality
  \[
    x= \sum_j \langle x, (\theta^*\theta)^{-1} (x_j) \rangle x_j =
       \sum_j \langle x,x_j \rangle (\theta^*\theta)^{-1} (x_j) \qquad
       {\it for} \quad x \in \mathcal H \, .
  \]
  In other words, the canonical bi-dual frame of a standard frame is the
  frame itself again. The frame transform $\theta'$ of the canonical dual
  frame $\{ (\theta^*\theta)^{-1} (x_j) : j \in \mathbb J \}$ equals
  $(\theta^*)^{-1}$, i.e.~the frame transform of the canonical dual frame
  maps $\mathcal H$ onto $\theta(\mathcal H) \subseteq l_2(A)$ acting like this
  operator.
\end{proposition}

\begin{proof}
By the definition of a canonical dual frame and by the results
of Theorem \ref{prop-canondual} above we have the equality
  \[
    x = \sum_j \langle x, (\theta^*\theta)^{-1} (x_j) \rangle x_j
  \]
for every $x \in \mathcal H$. Applying the invertible positive operator
$(\theta^*\theta)^{-1}$ to this equality we obtain the identity
  \begin{eqnarray*}
    (\theta^*\theta)^{-1} (x)
    & = & \sum_j \langle x, (\theta^*\theta)^{-1} (x_j) \rangle
            (\theta^*\theta)^{-1} (x_j) \\
    & = & \sum_j \langle (\theta^*\theta)^{-1}(x),x_j \rangle
            (\theta^*\theta)^{-1} (x_j)
  \end{eqnarray*}
for $x \in \mathcal H$. Since the operator $(\theta^*\theta)^{-1}$ is invertible
on $\mathcal H$ we can replace $(\theta^*\theta)^{-1}(x)$ by $x$, and the sought
equality turns out. By the uniqueness result of Theorem \ref{prop-canondual}
for the calculation of canonical dual frames and by the trivial equality
${\rm id}_{\mathcal H} = {\rm id}_{\mathcal H}^*{\rm id}_{\mathcal H}$ the
canonical bi-dual frame of a given standard frame equals the frame itself.
To calculate the frame transform $\theta'$ of the canonical dual frame
consider the special description of the identity map on $\mathcal H$
  \[
    x \stackrel{\theta'}{\longrightarrow}
    \{ \langle x,(\theta^*\theta)^{-1}(x_j) \rangle \}_{j \in \mathbb J}
    \stackrel{\theta^*}{\longrightarrow}
    \sum_j \langle x, (\theta^*\theta)^{-1}(x_j) \rangle x_j = x
  \]
($x \in \mathcal H$), cf.~Theorem \ref{th-reconstr}. Note, that $\{ \langle
x,(\theta^*\theta)^{-1}(x_j) \rangle \}_{j \in \mathbb J}$ belongs to
$P(l_2(A))$ since the operator $(\theta^*\theta)^{-1}$ is positive. The
equality shows $\theta'= (\theta^*)^{-1}$ as operators from $\mathcal H$ onto
$\theta(\mathcal H)$.
\end{proof}

The next proposition gives us the certainty that the relation between a frame
and its dual is symmetric. The equality tells us something about the relation
of the associated frame transforms. (Cf.~\cite[Prop.~1.13, 1.17]{HanLarson}.)

\begin{proposition}
  Let $\{ x_j : j \in \mathbb J \}$ and $\{ y_j : j \in \mathbb J \}$ be
  standard frames in a Hilbert $A$-module $\mathcal H$ with the property that
  they fulfil the equality $x = \sum_j \langle x,y_j \rangle x_j$ for every
  $x \in \mathcal H$.
  Then the equality $x=\sum_j \langle x,x_j \rangle y_j$ holds for every
  $x \in \mathcal H$, too.
  \newline
  Let $\theta_1$ and $\theta_2$ be the associated frame transforms of two
  frames $\{ x_j : j \in \mathbb J \}$ and $\{ y_j : j \in \mathbb J \}$ of
  $\mathcal H$, respectively. Then these two frames are duals to each other if
  and only if $\theta_2^*\theta_1= {\rm id}_{{\mathcal H}}$.
\end{proposition}

\begin{proof}
By Proposition \ref{prop-dualframe2} there exists a standard Riesz basis
$\{ f_j : j \in \mathbb J \}$ of a Hilbert $A$-module $\mathcal K$ and an
orthogonal projection $P$ such that $y_j=P(f_j)$ for $j \in {\mathbb J}$.
Since the sum $\sum_j \langle x,x_j \rangle \langle x_j,x \rangle$ is
norm-bounded we can define another adjointable operator $T: {\mathcal H} \to
{\mathcal K}$ by the formula $T(x) = \sum_j \langle x,x_j \rangle f_j$ for
$x \in \mathcal H$. Then $PT \in {\rm End}_A^*({\mathcal H})$ and $PT(x)=
\sum_j \langle x,x_j \rangle y_j$ for $x \in \mathcal H$. The following
equality holds for any $x \in \mathcal H$:
  \begin{eqnarray*}
    \langle x,x \rangle & = &
     \left\langle \sum_j \langle x,y_j \rangle x_j,x \right\rangle =
     \sum_n \langle x,y_j \rangle \langle x_j,x \rangle \\
     & = &
     \sum_j \langle x,x_j \rangle \langle y_j,x \rangle =
     \left\langle \sum_j \langle x,x_j \rangle y_j,x \right\rangle   \\
     & = &
     \langle PT(x),x \rangle
   \end{eqnarray*}
In the middle step we used the self-adjointness of $\langle x,x \rangle$.
As a result $PT$ is shown to be positive, and its square root to be an
isometry (cf.~\cite[Lemma 4.1]{Lance2}).
Since $PT = (PT)^{1/2}((PT)^{1/2})^*=((PT)^{1/2})^*(PT)^{1/2}=
{\rm id}_{{\mathcal H}}$ the operator $(PT)^{1/2}$ is at the same time a
unitary, and $PT={\rm id}_{{\mathcal H}}$. This demonstrates the first
assertion.

\smallskip
Now, let $x,y \in \mathcal H$, $\{ e_j : j \in \mathbb J \}$ be the standard
orthonormal Hilbert basis of $l_2(A)$ and $\{ x_j : j \in \mathbb J \}$ and
$\{ y_j : j \in \mathbb J \}$ be two frames of $\mathcal H$ with their
associated frame transforms $\theta_1,\theta_2$. We have the equality:
  \begin{eqnarray*}
    \langle \theta_1^*\theta_2 (x), y \rangle & = &
      \langle \theta_2(x),\theta_1(y) \rangle_{l_2(A)} =
      \left\langle \sum_j \langle x,y_j \rangle e_j, \sum_i \langle y,x_i
         \rangle e_i \right\rangle_{l_2(A)} \\
      & = &
      \sum_j \langle x,y_j \rangle \langle x_j,x \rangle =
      \left\langle \sum_j \langle x,y_j \rangle x_j , y \right\rangle \, .
  \end{eqnarray*}
Since $y \in \mathcal H$ is arbitrarily chosen the equality $\theta_1^*
\theta_2 (x) =\sum_j \langle x,y_j \rangle x_j$ turns out to hold for
every $x \in \mathcal H$. Therefore, $x=\sum_j \langle x,y_j \rangle x_j$ for
every $x \in \mathcal H$ if and only if $\theta_2^*\theta_1=
{\rm id}_{{\mathcal H}}$. With a reference to the definition of a dual frame
(see equation (\ref{eq-dualframe})) we are done.
\end{proof}

In contrast to the Hilbert space situation Riesz bases of Hilbert C*-modules
may possess infinitely many alternative dual frames due to the existence of
zero-divisors in the C*-algebra of coefficients, compare with
\cite[Cor.~2.26]{HanLarson}.

\begin{example}  \label{exam00} {\rm
   Let $A=l_\infty$ be the C*-algebra of all bounded complex-valued sequences
   and let ${\mathcal H}=c_0$ be the Hilbert $A$-module and two-sided ideal in
   $A$ of all sequences converging to zero. The $A$-valued inner product on
   $\mathcal H$ is that one inherited from $A$. Consider a maximal set of pairwise
   orthogonal minimal projections $\{ p_i : i \in {\mathbb Z} \}$ of $\mathcal H$.
   Since $x= \sum_i x p_i = \sum_i \langle x,p_i \rangle_A p_i$ for any $x \in
   \mathcal H$ and since the zero element admits a unique decomposition of this
   kind this set is a standard Riesz basis of $\mathcal H$, even an orthogonal
   Hilbert basis and a standard normalized tight frame at the same time.
   However, the Riesz basis $\{ p_i : i \in {\mathbb Z} \}$ possesses infinitely
   many alternate dual frames, for example $\{ p_i+p_{i+m} : i \in {\mathbb Z}
   \}$ for a fixed non-zero integer $m$. }
\end{example}

\begin{proposition}    \label{prop-optimal}
  Let $\{ x_j : j \in \mathbb J \}$ be a standard frame of a finitely or
  countably generated Hilbert $A$-module $\mathcal H$ over a unital C*-algebra
  $A$ that possesses more than one dual frame. Then for the canonical dual
  frame $\{ S(x_j) : j \in \mathbb J \}$ and for any other alternative dual
  frame $\{ y_j : j \in \mathbb J \}$ of the frame $\{ x_j : j \in \mathbb J
  \}$ the inequality
  \[
     \sum_{j} \langle x,S(x_j) \rangle\langle S(x_j),x \rangle
     \leq \sum_{j} \langle x,y_j \rangle \langle y_j,x \rangle
  \]
  is valid for every $x \in \mathcal H$. Beside this, equality holds precisely
  if $S(x_j)=y_j$ for every $j \in \mathbb J$.

  More generally, whenever $x = \sum_{j \in \mathbb J} a_jx_j$ for certain
  elements $a_j \in A$ and $\sum_{j \in \mathbb J} a_ja_j^*$ is bounded in
  norm we have
  \[
     \sum_{j} a_ja_j^* =
     \sum_{j} \langle x,S(x_j) \rangle\langle S(x_j),x \rangle +
     \sum_{j} (a_j -\langle x,S(x_j) \rangle)(a_j -\langle x,S(x_j) \rangle)^*
  \]
  with equality in case $a_j = \langle x,S(x_j) \rangle$ for every $j \in
  \mathbb J$. Moreover, the minimal value of the summands $a_j^*a_j$ is
  admitted for $a_j=\langle x,S(x_j)\rangle$ for each $j \in \mathbb J$
  separately. (Cf.~Example \ref{exam00}.)
\end{proposition}

\begin{proof}
We begin with the proof of the first statement. The convergence of the sums
in the inequality above follows from the properties of the frame transforms
and of the frame operators.
If the standard frames $\{ S(x_j) : j \in \mathbb J \}$ and $\{ y_j : j \in
\mathbb J \}$ are both dual frames of $\{ x_j : j \in \mathbb J \}$ then the
equalities
  \[
    x= \sum_{j} \langle x,S(x_j) \rangle x_j
     = \sum_{j} \langle x,y_j \rangle x_j
  \]
are valid for every $x \in \mathcal H$. Subtracting one sum from the other,
applying the operator $S$ to the result and taking the $A$-valued inner
product with $x$ from the right we obtain
  \[
    0 = \sum_{j} \langle x,y_j-S(x_j) \rangle \langle S(x_j),x
    \rangle
  \]
for every $x \in \mathcal H$. Therefore,
  \begin{eqnarray*}
    \sum_{j} \langle x,y_j \rangle \langle y_j,x \rangle
    & = & \sum_{j} \langle x,y_j -S(x_j)+S(x_j) \rangle \langle
          y_j -S(x_j)+S(x_j),x \rangle  \\
    & = & \sum_{j} \langle x,S(x_j) \rangle \langle S(x_j),x
       \rangle + \sum_{j} \langle x,y_j-S(x_j) \rangle \langle
       y_j-S(x_j),x \rangle   \, ,
  \end{eqnarray*}
demonstrating the first part of the stressed for assertion since every
summand is a positive element of $A$.

Now suppose $x \in \mathcal H$ has two decompositions
$x = \sum_j \langle x,S(x_j) \rangle x_j = \sum a_j x_j$ with coefficients
$\{ a_j \}_j \in l_2(A)$, where the index set $\mathbb J$ has to be identified
with $\mathbb N$ to circumvent extra discussions about conditional and unconditional
convergence of series. Then the equality
\[
   0 = \sum_j (\langle x,S(x_j) \rangle - a_j) \langle x_j,S(x) \rangle
     = \sum_j (\langle x,S(x_j) \rangle - a_j) \langle S(x_j),x \rangle
\]
holds by the self-adjointness of $S$. Consequently,
\begin{eqnarray*}
  \langle \{ a_j \}_j, \{ a_j \}_j \rangle_{l_2(A)}  & = &
    \langle \{ \langle x,S(x_j) \rangle \}_j , \{ \langle x,S(x_j) \rangle \}_j
    \rangle_{l_2(A)} + \\
    & \,\, + &
    \langle \{ \langle x,S(x_j) \rangle - a_j \}_j , \{ \langle x,S(x_j) \rangle
    - a_j \}_j \rangle_{l_2(A)} \, ,
\end{eqnarray*}
and by the positivity of the summands the minimal value of $a_j a_j^*$ is
admitted for
   \linebreak[4]
$a_j=\langle x,S(x_j) \rangle$ for each $j \in {\mathbb J}$ separately.
\end{proof}

The optimality principle allows to investigate the stability of the frame
property to be standard under changes of the $A$-valued inner product on
Hilbert C*-modules. The result is important since countably generated
Hilbert C*-modules may possess non-adjointable bounded module isomorphisms
and, consequently, a much wider variety of C*-valued inner products inducing
equivalent norms to the given one than Hilbert spaces use to admit,
cf.~\cite{Frank:93}.

\begin{corollary}
   Let $A$ be a unital C*-algebra, $\mathcal H$ be a finitely or countably
   generated Hilbert $A$-module with $A$-valued inner product $\langle .,.
   \rangle_1$ and $\{ x_j : j \in \mathbb J \} \subset \mathcal H$ be a
   standard frame. Then $\{ x_j : j \in \mathbb J \}$ is a standard frame
   with respect to another $A$-valued inner product $\langle .,. \rangle_2$
   on $\mathcal H$ that induces an equivalent norm to the given one, if and
   only if there exists an adjointable invertible bounded operator $T$ on
   $\mathcal H$ such that $\langle .,. \rangle_1 \equiv \langle T(.),T(.)
   \rangle_2$. In that situation the frame operator $S_2$ of $\{ x_j : j
   \in \mathbb J \}$ with respect to $\langle .,. \rangle_2$ commutes with the
   inverse of the frame operator $S_1$ of $\{ x_j : j \in \mathbb J \}$ with
   respect to $\langle .,. \rangle_1$.
\end{corollary}

\begin{proof}
Suppose the frame $\{ x_j \}_{j \in \mathbb J}$ is standard with respect to
both the inner products on $\mathcal H$. For $x \in \mathcal H$ we have two
reconstruction formulae $x = \sum_j \langle x,S_1(x_j) \rangle_1 x_j$ and
$x = \sum_j \langle x,S_2(x_j) \rangle_2 x_j$. By the optimality principle we
obtain the equality $\langle S_1(x),x_j \rangle_1 = \langle x,S_1(x_j)
\rangle_1 = \langle x,S_2(x_j) \rangle_2 = \langle S_2(x),x_j \rangle_2$
that is satisfied for any $x \in \mathcal H$ and $j \in \mathbb J$, see
Proposition \ref{prop-optimal}. Let $y \in \mathcal H$. Multiplying by
$\langle S_1(x_j),y \rangle_1$ from the right and summing up over $j \in
\mathbb J$ we arrive at the equality $\langle S_1(x),y \rangle_1 = \langle
S_2(x),y \rangle_2$ that has to be valid for any $x,y \in \mathcal H$.
Therefore, $0 \leq \langle z,z \rangle_1 = \langle z,(S_2 S_1^{-1})(z)
\rangle_2$ for any $z \in \mathcal H$ forcing $(S_2S_1^{-1})$ to be
self-adjoint and positive by \cite[Lemma 4.1]{Lance2}. In particular, the
operators commute since $S_2$ itself is positive with respect to the second
inner product by construction.
So we can take the square root of this operator in the C*-algebra of all
adjointable bounded module operators on $\mathcal H$ as the particular operator
$T$ that relates the $A$-valued inner products one to another by
$\langle .,. \rangle_1 \equiv \langle T(.),T(.) \rangle_2$.

Conversely, if both the $A$-valued inner products on $\mathcal H$ are
related as $\langle .,. \rangle_1 \equiv \langle T(.),T(.) \rangle_2$ for
some adjointable invertible bounded operator $T$ on $\mathcal H$ then the
frame operators fulfil the equality $S_1=T^*S_2T$, and the frame
$\{ x_j : j \in \mathbb J \}$ is standard with respect to both the inner
products.
\end{proof}

Different alternate duals of a standard frame cannot be similar
or unitarily equivalent in any situation, so we reproduce a Hilbert space
result (\cite[Prop.~1.14]{HanLarson}).

\begin{proposition}     \label{prop-6.7}
  Suppose, for a given standard frame $\{ x_j : j \in \mathbb J \}$ of a
  Hilbert $A$-module $\mathcal H$ over a unital C*-algebra $A$ there exist
  two standard alternate dual frames $\{ y_j : j \in \mathbb J \}$ and
  $\{ z_j : j \in \mathbb J \}$ which are connected by an invertible adjointable
  operator $T \in {\rm End}_A({\mathcal H})$ via $z_j=T(y_j)$, $j \in
  {\mathbb J}$. Then $T={\rm id}_{{\mathcal H}}$. \newline
  In other words, two different standard alternate dual frames of a given
  frame are not similar or unitarily equivalent.
\end{proposition}

\begin{proof}
Suppose $z_j=T(y_j)$ for $j \in \mathbb J$ and an adjointable invertible
operator $T$. Let us count the values of the adjoint operator $T^*$ of $T$.
We have $T^*(x) = \sum_j \langle T^*(x),y_j \rangle x_j = \sum_j \langle
x,T(y_j) \rangle x_j = x$ for every $x \in {\mathcal H}$ by the dual frame
property. Consequently, $T$ equals the identity operator on $\mathcal H$.
\end{proof}

We conjecture that the restriction to $T$ to be adjointable may be dropped
preserving the assertion of the proposition. To check this techniques described
at the appendix of the present paper might be helpful.

\medskip
Since for every orthogonal projection $P$ on a Hilbert C*-module $\mathcal H$
and every standard frame $\{ x_j : j \in \mathbb J \}$ of $\mathcal H$ the
sequence $\{ P(x_j) : j \in \mathbb J \}$ is a standard frame for the Hilbert
C*-submodule $P({\mathcal H})$ the natural question is whether the canonical
dual frame of this frame $\{ P(x_j) : j \in \mathbb J \}$ would be equal
to the projected canonical dual frame of $\{ x_j : j \in \mathbb J \}$, or not.
If the frame $\{ x_j : j \in \mathbb J \}$ is tight then we get a global
affirmative answer. However, if $\{ x_j : j \in \mathbb J \}$ is not tight
then the projection $P$ has to commute with the frame operator $S$ related to
the frame $\{ x_n \}$ to guarantee the square of these mappings to commute.
However, every orthogonal projection of the canonical dual frame is still a
standard alternate dual frame because
  \[
    x = P(x) = \sum_j \langle P(x),S(x_j) \rangle P(x_j)
             = \sum_j \langle x,PS(x_j) \rangle P(x_j)
  \]
for every $x \in P({\mathcal H})$. Unfortunately, the set of orthogonal
projections on a Hilbert C*-module may be very small, in extreme cases
reducing to the zero and the identity ope\-rator. Nevertheless, for existing
projection operators the analogous to \cite[Prop.~1.15, Cor.~1.16]{HanLarson}
facts hold:

\begin{proposition}
  Let $\{ x_j : j \in \mathbb J \}$ be a standard frame of a Hilbert C*-module
  $\mathcal H$ and $S_x \geq 0$ be its frame operator. If $P$ is an orthogonal
  projection on $\mathcal H$ then the frame operator of the projected frame
  $\{ P(x_j) : j \in \mathbb J \}$ is $S_{P(x)}= PS_x$ if and only if
  $PS_x=S_xP$.
  \newline
  The standard frame $\{ x_j : j \in \mathbb J \}$ is tight if and only if
  $S_x$ equals the identity operator multiplied by the inverse of the frame
  bound. In this situation the equality $S_{P(x)}= PS_x$ is fulfilled for
  every orthogonal projection on $\mathcal H$. Conversely, the latter
  condition alone does not imply the frame $\{ x_j : j \in \mathbb J \}$ to
  be tight, in general.
\end{proposition}

\begin{proof}
Considering the first pair of equivalent conditions the product
of the two positive elements $S_x$ and $P$ of the C*-algebra ${\rm End}_A^*
({\mathcal H})$ can only be positive if they commute. Consequently, $S_{P(x)}=PS_x$
forces them to commute since $S_{P(x)} \geq 0$ by construction,
cf.~Theorem \ref{prop-canondual}.

Conversely, if $PS_x=S_xP$ then by the equality $x = \sum_n \langle x,S_x(x_n)
\rangle x_n$ for $x \in \mathcal H$ we obtain
  \begin{eqnarray*}
    P(x) & = & \sum_n \langle P(x),S_x(x_n) \rangle P(x_n)
           =   \sum_n \langle P(x),PS_x(x_n) \rangle P(x_n) \\
         & = & \sum_n \langle P(x),S_xP(x_n) \rangle P(x_n)
           =   \sum_n \langle P(x),(S_xP)(P(x_n)) \rangle P(x_n)  \, .
  \end{eqnarray*}
By the positivity of $PS_x=S_xP$, by the free choice of $x \in \mathcal H$,
by Theorem \ref{prop-canondual} and by Proposition \ref{prop-6.7} the equality
$S_{P(x)}=PS_x$ turns out to hold.

The second statement is nearly obvious. Since there are C*-algebras with
very small sets of projections, like $A={\rm C}([0,1])$, the property of the
frame operator $S_x$ of an one-element frame $\{ x=a \} \in A$ to commute with
any projection $P \in {\rm End}_A^*(A)$ does certainly not imply the frame
to be tight. In our example any invertible element $a \in A$ has this
property despite of its possibly unequal to one norm or frame bounds.
\end{proof}

We add a few more remarks on the properties of the frame
transform $\theta$ and of the operator $(\theta^*)^{-1}: \mathcal H \to
\theta(\mathcal H)$. For this aim consider the operator $R=\theta S$.
This operator $R$ has the property that $R^*\theta = {\rm id}_{\mathcal H}
= \theta^*R$ by the definition of $S$ and $\theta$, cf.~Theorem
\ref{prop-canondual}. Moreover, the equality $\theta(R^*\theta) = (\theta R^*)
\theta = \theta$ and the injectivity of $\theta$ imply $\theta R^*=P$ on
$l_2(A)$. Also, $\theta R^*=R \theta^*$ as can be easily verified. Therefore,
$R^*$ restricted to $\theta(\mathcal H)$ is an inverse to the operator $\theta$,
and $R$ is an inverse of the operator $\theta^*$ if $\theta^*$ has been
restricted to $\theta(\mathcal H)$. So, alternative descriptions of the
situation can be given in terms of a quasi-inverse operator for the extension
of the frame transform $\theta$ to an operator on $\mathcal H \oplus l_2(A)$.
%%%%%%%%%%%%%%%%%%%%%%%%%%%%%%%%%%%%%%%%%%%%%%%%%%%%%%%%%%%%%%%%%%%%%%%%%%%%%%%

\section{A classification result}

We would like to get a better understanding of the unitary and similarity
equivalence classes of frames in a Hilbert C*-module with orthogonal basis.
Comparing the result with the results of section \ref{sec-complement} we get
general insights into necessary conditions for frame equivalence in Hilbert
C*-modules, even in the absence of an orthogonal Hilbert basis for them.
For the Hilbert space situation we refer to \cite[Prop.~2.6]{HanLarson}.

\begin{proposition}  \label{prop-equiv-id}      %Prop. 2.6 in {HL}, better.
  Let $A$ be a C*-algebra and $\mathcal H$ be a countably generated Hilbert
  $A$-module with orthogonal Hilbert basis $\{ f_j : j \in {\mathbb J} \}$.
  For two orthogonal projections $P,Q \in {\rm End}_A^*({\mathcal H})$ set
  ${\mathcal M} = P({\mathcal H})$ and ${\mathcal N} = Q({\mathcal H})$. Let
  the sequences $\{ x_j=P(f_j) : j \in {\mathbb J} \}$ and $\{ y_j=Q(f_j) :
  j \in {\mathbb J} \}$ be the derived standard normalized tight frames for
  $\mathcal M$ and $\mathcal N$, respectively. Then the frames $\{ x_j : j
  \in \mathbb J \}$ and $\{ y_j : j \in \mathbb J \}$ are similar if and only
  if they are unitarily equivalent, if and only if $P=Q$ and the frames
  coincide elementwise.
\end{proposition}

\begin{proof}
Suppose, there exists an adjointable invertible bounded $A$-linear operator
\linebreak[4]
$T: {\mathcal M} \to {\mathcal N}$ with $T(x_j)=y_j$ for every
$j \in {\mathbb J}$. Continuing the operator $T$ and its adjoint on the
orthogonal complements of $\mathcal M$ and $\mathcal N$, respectively, as the
zero operator we obtain an adjointable bounded $A$-linear operator $T$ on
$\mathcal H$ that possesses a polar decomposition in ${\rm End}_A^*
({\mathcal H})$, $T=V \cdot |T|$ (cf.~\cite[Th.~15.3.7]{Wegge-Olsen}). The
partial isometry $V$ has the property $VV^*=Q$, $V^*V=P$ by construction.
Furthermore, $y_j = T(x_j) = V \cdot |T| (x_j)$. Since $\{ y_j : j \in \mathbb
J \}$ is normalized tight, since $V$ is an isometry of ${\mathcal M}$ with
${\mathcal N}$ and because $T$ is invertible the standard frame $\{ |T|(x_j)
: j \in \mathbb J \}$ has to be a standard normalized tight frame for
$\mathcal M$. Also, $|T|= {\rm id}_{{\mathcal M}}$. So similarity implies
unitary equivalence.

Let us continue with the partial isometry $V$ obtained above. The operator $V$
canonically arises if we suppose the frames $\{ x_j : j \in \mathbb J \}$ and
$\{ y_j : j \in \mathbb J \}$ to be unitarily equivalent. Since $V = VP$ we
obtain $V(f_j) = VP(f_j) =Q(f_j)$ for every $j \in {\mathbb J}$. Since
$\{ f_j : j \in {\mathbb J} \}$ is an (orthogonal) Hilbert basis of $\mathcal
H$ we find $V=Q$ and hence, $P=Q$ and $x_j=y_j$ for every $j \in {\mathbb J}$.
\end{proof}

The next theorem and the derived from it corollary give us a criterion on
similarity and unitary equivalence of frames in Hilbert C*-modules. They
generalize \cite[Cor.~2.8, 2.7]{HanLarson} and \cite[Th.~B]{Holub:97} and tie
these observations together.

\begin{theorem} \label{th-equiv-frametransform}
  Let $A$ be a unital C*-algebra and $\{ x_j : j \in {\mathbb J} \}$ and $\{ y_j :
  j \in {\mathbb J} \}$ be standard normalized tight frames of Hilbert $A$-modules
  ${\mathcal H}_1$ and ${\mathcal H}_2$, respectively. Then the frames $\{ x_j
  : j \in \mathbb J \}$ and $\{ y_j : j \in \mathbb J \}$ are unitarily
  equivalent if and only if their frame transforms $\theta_1$ and $\theta_2$
  have the same range in $l_2(A)$, if and only if the sums $\sum_j a_j x_j$ and
  $\sum_j a_j y_j$ equal zero for exactly the same Banach A-submodule of
  sequences $\{ a_j : j \in {\mathbb J} \}$ of $l_2(A)'$.
  \newline
  Similarly, two standard frames $\{ x_j : j \in {\mathbb J} \}$ and $\{ y_j :
  j \in {\mathbb J} \}$ of Hilbert $A$-modules ${\mathcal H}_1$ and ${\mathcal H}_2$,
  respectively, are similar if and only if their frame transforms $\theta_1$
  and $\theta_2$ have the same range in $l_2(A)$, if and only if
  the sums $\sum_j a_j x_j$ and $\sum_j a_j y_j$ equal zero for exactly the same
  Banach A-submodule of sequences $\{ a_j : j \in {\mathbb J} \}$ of $l_2(A)'$.
\end{theorem}

\begin{proof}
If we assume that the frame transforms $\theta_1$, $\theta_2$ corresponding
to the two initial standard normalized tight frames have the same range in
$l_2(A)$ then the orthonormal projection $P$ of $l_2(A)$ onto this range
$\theta_1({\mathcal H}_1) \equiv \theta_2({\mathcal H}_2)$ maps the elements
of the standard orthonormal basis $\{ e_j : j \in \mathbb J \}$ of $l_2(A)$
to both $\theta_1(x_j)=\theta_2(y_j)$, $j \in \mathbb J$, by the construction
of a frame transform, cf.~Proposition \ref{prop-complement} and Theorem
\ref{th-reconstr}. Then
\[
\langle x_j,x_j \rangle_1 = \langle \theta_1(x_j),\theta_1(x_j) \rangle_{l_2}
= \langle \theta_2(y_j),\theta_2(y_j)\rangle_{l_2} = \langle y_j,y_j \rangle_2
\]
for every $j \in {\mathbb J}$, and the mapping $U: {\mathcal H}_1 \to {\mathcal H}_2$
induced by the formula $U(x_j)=y_j$ for $j \in {\mathbb J}$ is a unitary
isomorphism since the sets $\{ x_j \}$ and $\{ y_j \}$ are sets of generators
of ${\mathcal H}_1$ and ${\mathcal H}_2$, respectively. Moreover, the set of bounded
$A$-linear functionals on $l_2(A)$ annihilating the ranges of the frame
transforms $\theta_1$, $\theta_2$ are exactly the same and can be identified
with a Banach $A$-submodule of $l_2(A)'$.

The converse statement for standard normalized tight frames follows directly
from Proposition \ref{prop-equiv-id}.

\smallskip
If we suppose merely $\{ x_j : j \in \mathbb J \}$ and $\{ y_j : j \in \mathbb
J \}$ to be standard frames in ${\mathcal H}_1$ and ${\mathcal H}_2$,
respectively, then the assumption $\theta_1({\mathcal H}_1) \equiv \theta_2
({\mathcal H}_2)$ yields $P(e_j)=\theta_1(x_j)=\theta_2(y_j)$ again,
cf.~Theorems \ref{th-reconstr} and \ref{th-orthogonal}. Consequently, there is
an adjointable invertible bounded operator $T \in {\rm End}_A({\mathcal H}_1,
{\mathcal H}_2)$ with $T(x_j)=y_j$ for $j \in {\mathbb J}$ by the injectivity
of frame transforms.
\end{proof}

\begin{corollary}
  Let $A$ be a unital C*-algebra.
  Let ${\mathbb J}$ be a countable (or finite, resp.) index set. The set of
  unitary equivalence classes of all standard normalized tight frames indexed
  by ${\mathbb J}$ is in one-to-one correspondence with the set of all
  orthogonal projections on the Hilbert $A$-module $l_2(A)$ (or $A^{|{\mathbb
  J}|}$, resp.).
  Analogously, the set of similarity equivalence classes of all frames
  indexed by ${\mathbb J}$ is in one-to-one correspondence with the set of all
  orthogonal projections on the Hilbert $A$-module $l_2(A)$
  (or $A^{|{\mathbb J}|}$, resp.). The one-to-one correspondence can be
  arranged fixing an orthonormal Hilbert basis of $A^{{\mathbb J}|}$ or
  $l_2(A)$, respectively.
\end{corollary}

The established interrelation allows to transfer the partial order structure
of projections as well as homotopy and other topological properties of the
set of projections to properties of equivalence classes of standard frames.
The resulting structures may strongly depend on the choice of some
orthonormal Hilbert basis realizing the correspondence. However, the partial
order does not depend on the choice of the orthonormal Hilbert basis since
orthonormal Hilbert bases of $l_2(A)$ (or of $A^{|{\mathbb J}|}$) are unitarily
equivalent. Despite the special situation for Hilbert spaces $\mathcal H$ the
C*-algebra ${\rm End}_A^*(l_2(A))$ has a partial ordered subset of projections
which lacks the lattice property for many C*-algebras $A$.

%%%%%%%%%%%%%%%%%%%%%%%%%%%%%%%%%%%%%%%%%%%%%%%%%%%%%%%%%%%%%%%%%%%%%%%%%%%%%%%
\section{Final remarks}

\smallskip
We would like to add some remarks on non-standard frames in C*-algebras and
Hilbert C*-modules. As we already mentioned in the introduction a good theory
can be developed for non-standard frames in self-dual Hilbert C*-modules over
von Neumann algebras or monotone complete C*-algebras since a well-defined
concept of a generalized Hilbert basis exists for that class of Hilbert
C*-modules, cf.~\cite{Pa1,BDH,Fr7,DH,Fr95}.
However, because of numerous Hilbert C*-module isomorphisms in this class
non-trivial examples may be only obtained, first, in the case of finite
W*-algebras of coefficients or secondly, for cardinalities of the index set
of the frame greater than the cardinality of every decomposition of the
identity into a sum of pairwise orthogonal and equivalent to one projections
in the complementary case of infinite W*-algebras of coefficients. The target
space for the frame transform is always $l_2(A,I)'$ for an index set $I$ of
the same cardinality as the index set $\mathbb J$ of the frame under
consideration. The first steps towards a frame theory for self-daual Hilbert
W*-modules can be found in a paper by Y.~Denizeau and J.-F.~Havet \cite{DH}
where a weak reconstruction formula appears.

In case of non-standard frames for Hilbert C*-modules over general C*-algebras
$A$ we have the difficulty to define a proper target space for the frame
transform where the image of the frame transform becomes a direct summand. The
choice of the C*-dual Banach $A$-module $l_2(A,I)'$ for a suitable index set
$I$ of the same cardinality as the index set of the frame may not always be
the right choice since the C*-dual Banach $A$-module of the initial Hilbert
C*-module carrying the frame set may not fit into $l_2(A,I)'$. The latter
phenomenon is mainly caused by the sometimes complicated multiplier theory
of ideals of $A$. A better candidate for the target space seems to be
$l_2(A^{**},I)'$ where $A^{**}$ denotes the bidual von Neumann algebra of $A$.
To embed the original Hilbert $A$-module ${\mathcal H}$ into $l_2(A^{**},I)'$
by a frame transform we have to enlarge ${\mathcal H}$ to an Hilbert
$A^{**}$-module by the techniques described in the appendix and afterwards
to 'self-dualize' it as described in \cite{Pa1}. The frame will preserve its
properties, i.e.~the frame will still be a frame for the larger Hilbert
$A^{**}$-module with the same frame bounds. For tight frames we obtain a
proper reconstruction formula with weak convergence of the occuring sum that
can be restricted to the original module ${\mathcal H}$ in such a way that any
trace of the made extensions vanishes. In particular, non-standard tight frames
are always generator sets in a weak sense. However, the frame transform is
only a modified one and does not map ${\mathcal H}$ to a direct summand of
$l_2(A^{**},I)'$.
(An alternative view on these facts can be given using linking C*-algebra
techniques.)

To make use of the complete boundedness of bounded C*-module maps between
Hilbert C*-modules and of injectivity properties of objects one could also
consider to take the atomic part of $A^{**}$ or the injective envelope $I(A)$
of $A$ instead of $A^{**}$ and to repeat the construction presented in the
appendix appropriately. This would let to operator space and operator
module methods. All in all we can say that a general theory of non-standard
frames in Hilbert C*-modules and C*-algebras doesn't exist at present. Steps
towards such a theory have to involve results from Banach space and
operator space theory, as well as from operator and operator algebra theory.

\begin{problem}  {\rm
   Whether every Hilbert C*-module over a unital C*-algebra admits a normalized
   tight frame, or not?  }
\end{problem}

%%%%%%%%%%%%%%%%%%%%%%%%%%%%%%%%%%%%%%%%%%%%%%%%%%%%%%%%%%%%%%%%%%%%%%%%%%%%%%%
\section{Appendix}

\smallskip
In proofs we need a canonical construction for a canonical switch from a given
Hilbert $A$-module $\mathcal M$ to a bigger Hilbert $A^{**}$-module ${\mathcal M}^\#$
while preserving many useful properties and guaranteeing the existence and
uniqueness of extended operators and $A$-($A^{**}$-)valued inner products.
The much better properties of Hilbert W*-modules in comparison to general
Hilbert C*-modules (cf.~\cite{Pa1}) and facts from non-commutative
topology form the background for such a manner of changing objects.

\begin{remark}  \label{***} {\rm  (cf.~H.~Lin \cite[Def.~1.3]{Lin:90/2},
                             \cite[\S 4]{Pa1})   \newline
    Let $\{ {\mathcal M},\langle .,. \rangle \}$ be a left pre-Hilbert
    $A$-module over a fixed C*-algebra $A$. The algebraic tensor
    product $A^{**} \odot {\mathcal M}$ becomes a left
    $A^{**}$-module defining the action of $A^{**}$ on its
    elementary tensors by the formula $a(b \otimes h) = a b \otimes h$ for
    $a,b \in A^{**}$, $h \in \mathcal M$. Setting
    \[
    \left[ \sum_i a_i \otimes h_i , \sum_j b_j \otimes g_j \right] =
    \sum_{i,j} a_i \langle h_i, g_j \rangle b_j^*
    \]
    on finite sums of elementary tensors we obtain a degenerate
    $A^{**}$-valued inner pre-product. Factorizing $A^{**}
    \odot {\mathcal M}$ by $N= \{ z \in A^{**} \odot {\mathcal M} : [z,z]=0
    \}$ we obtain a pre-Hilbert $A^{**}$-module subsequently denoted by
    ${\mathcal M}^{\#}$. The pre-Hilbert $A^{**}$-module ${\mathcal M}^{\#}$
    contains $\mathcal M$ as a $A$-submodule.
    If $\mathcal M$ is Hilbert, then ${\mathcal M}^{\#}$ is Hilbert, and vice versa.
    The transfer of self-duality is more difficult. If $\mathcal M$ is
    self-dual, then ${\mathcal M}^{\#}$ is also self-dual by
    \cite[Th.~6.4]{Frank:93} and \cite{Pa1,Fr1}.

    \begin{problem} {\rm
    Suppose, the underlying C*-algebra $A$ is unital.
    Does the property of ${\mathcal M}^{\#}$ of being self-dual imply that
    $\mathcal M$ was already self-dual?}
    \end{problem}

    A bounded $A$-linear operator $T$ on $\mathcal M$ has a unique extension
    to a bounded $A^{**}$-linear operator on ${\mathcal M}^{\#}$ preserving
    the operator norm, (cf.~\cite[Def.~1.3]{Lin:90/2}).   }
\end{remark}

%%%%%%%%%%%%%%%%%%%%%%%%%%%

%%%%%%%%%%%%%%%%%%%%%%%%%%%

\end{document}